\documentclass{article}
\usepackage{amsmath,amssymb,graphicx,stmaryrd,enumerate}

\newtheorem{theorem}{Theorem}
\newtheorem{definition}{Definition}
\newtheorem{proposition}{Proposition}
\newtheorem{corollary}{Corollary}
\newtheorem{lemma}{Lemma}
\newtheorem{remark}{Remark}

\newcommand{\rI}{\rm I}
\newcommand{\rII}{\rm II}
\newcommand{\rIII}{\rm III}
\newcommand{\st}{\textnormal{ s.t. }}
\newcommand{\bl}{\textnormal{ }}
\newcommand{\noask}{\noalign{\smallskip}}
\newcommand{\eop}{\hfill $\Box$\medskip}

\numberwithin{equation}{section}

\begin{document}
\title{On the Generating Functionals of a Class
of Random Packing Point Processes}
\author{Tien Viet Nguyen\footnote{Qualcomm, USA}\quad and \quad 
\hspace{-.2cm}Fran\c{c}ois Baccelli\footnote{UT Austin, USA and INRIA, France}}
\maketitle
\begin{abstract}
Consider a symmetrical conflict relationship between the points of a point
process. The Mat\'ern type constructions provide a generic way of
selecting a subset of this point process which is conflict-free.
The simplest one consists in keeping only conflict-free points.
There is however a wide class of Mat\'ern type processes based on
more elaborate selection rules and providing larger sets
of selected points. The general idea being that if a point is
discarded because of a given conflict, there is no need to discard 
other points with which it is also in conflict. The ultimate selection
rule within this class is the so called Random Sequential Adsorption,
where the cardinality of the sequence of conflicts allowing
one to decide whether a given point is selected is {\em not} bounded. 
The present paper provides a sufficient condition on the span
of the conflict relationship under which all the above point processes
are well defined when the initial point process is Poisson.
It then establishes, still in the Poisson case, a set of
differential equations satisfied by
the probability generating functionals of these Mat\'ern type
point processes.
Integral equations are also given for the
Palm distributions.
\end{abstract}

\section{Introduction}\label{S:Intro}

This paper is concerned with systems $\phi$ made of 
a collection of points in $\mathbb{R}^{d}$,
where points are equipped with a pairwise {\em conflict relation}.
As an example, one can think of the locations of dominant male lions
in a region.
These lions compete with each other for resources, 
so that they do not like to be too close to each other. Hence, it makes sense
to say that there is a conflict between two points representing the locations
of two such lions if the distance between them is, say, smaller
some distance $d$.

These questions play an important role in physics, chemistry,
material science, etc., where they have been used to describe systems 
with hard-core interactions such as reactions on polymer chains \cite{Gon74},
chemisorption on a single-crystal surface  \cite{Evan93} and adsorption in colloidal \cite{Senger00}
systems. In these problems, each object in the system occupies some space,
and two points with overlapping occupied space have a conflict.
Problems of this type also arise in forestry \cite{Matern60} and in wireless communications \cite{Busson08}.

The paper is centered on what happens after
conflicts are resolved (e.g. once, for all pairs of  
conflicting lions, one has eliminated the other).
More precisely, we focus our attention on the
three classical constructions of {\em conflict-free subsets}
$\phi_{j}$ ($j=$ I, II, III) of $\phi$ which were proposed by 
Mat\'ern \cite{Matern86} and which are described below.

The setting will be that where the locations
of points and the conflict relations are random. The conflict system
will hence be a random point process (PP) equipped with a random
conflict relation. The Mat\'ern constructions then lead to a family
of PPs $\Phi_{j}$ ($j=1,2,\ldots,\infty$), which are
almost surely (a.s.) conflict-free
when $j$ is either even or equal to $\infty$. 

The paper is structured as follows.
The models are introduced in Section \ref{S:Model}.
The existence issue is discussed in Section \ref{S:type III extend}.
The dynamical view point on these processes, which
is instrumental to establish the differential equations, 
is introduced in Section \ref{secprain}.
The differential equations satisfied by the
generating functionals are established in Section
\ref{S:Matern pgfl}, and in Section \ref{S: Palm int} for the Palm case.

\section{State of the Art and Contributions}
In dimension 1, the first model of this kind is the
car parking model which was independently studied by
A. R\'enyi \cite{Renyi58} and by
H. Dvoretzky and A. Robbins \cite{Robbins85}.
In this model, cars of fixed length are parked in
the same manner as in the $\infty$-Mat\'ern model (see below).
Consider an observation window $[0,x]$ and let
$N(x)$ be the number of cars parked in this window when there
is an infinite number of cars to be parked (saturated regime).
A. R\'enyi showed that $N(x)$ satisfies the law of large number (LLN):
\[
\lim_{x \rightarrow \infty}\frac{N(x)}{x}=C \approx 0.74759 \mbox{ a.s. },
\]
where $C$ is called the packing density.
H. Dvoretzky and A. Robbins \cite{Robbins85}
sharpened this result to a central limit theorem (CLT):
\[
\frac{N(x)-Cx}{Var(N(x))}\rightarrow \mathcal{N}(0,1) \mbox{ in distribution as } x\rightarrow \infty. 
\]
Various extensions of the above models were considered like the non-saturated
regime (the number of cars to be parked is finite), random car lengths, etc.
\cite{Mullooly68,Cof98}.
The latter is also known under the name random interval packing and
has many applications in resource allocation in communication theory.
For the above models, the obtained results concern the packing
density, the LLN, the CLT, the distribution of packed intervals
and that of vacant intervals.

For dimension more than $1$, the most noticeable advance
in this field is a series of papers by M.D. Penrose,
J.E. Yukich and Y. Baryshnikov. Based on a general LLN and CLT for
\emph{stabilizing} functionals, the LLN and CLT were established
for the $n$ dimensional RSA model in the non-saturated \cite{Pen02} regime.
Y. Baryshnikov and J. E. Yukich
\cite{Yukich03} later strengthened the above results by proving that,
in the thermodynamic limit, the spatial distribution of the p.p.
induced by the RSA model converges to that of a Gaussian field
after a suitable rescaling in the non-saturated regime.
The LLN and CLT for the $n$ dimensional RSA model in the saturated 
regime was proven by T. Schreiber \emph{et al} \cite{Schreiber07}.

In general, characterizing the distribution
of the point processes induced by these models
beyond first and second moment measures is still an open problem.
The results in this paper shed some light on this characterization.
The contributions are twofold:
we first give conditions under which each $\Phi_{j}$
is well defined; under these conditions, we then characterize
its probability generating functional (p.g.fl)
as the unique solution of a system of
functional differential equations. Using the
relation between the Palm and non-Palm distribution of a PP,
we also show that the p.g.fl of $\Phi_{j}$ under its Palm
distributions satisfies a related system of integral equations.
To the best of our knowledge, these systems of equations 
are new.

\section{Mat\'ern's Constructions and Extensions}\label{S:Model}
Let $c$ be a symmetric, non-reflexive, $\{0,1\}$--valued relation between the points of $\phi$. For any $x$ and $y$ in $\phi$, $c(x,y)=1$ means that $x$ and $y$ conflict with each other .

 \sl{The Mat\'ern type I construction} \rm builds the conflict-free system $\phi_{\rI}$ by removing from $\phi$ all objects which conflict with at least another object in $\phi$. In other words,
\begin{align}
 \phi_{\rI}=\{x \in \phi \st \textnormal{for all other } y \in \phi, \bl c(x,y)=0 \}.
\end{align}
 For convenience, we refer to a construction of the last type as a \emph{conflict resolution mechanism}. It is important to bear in mind that this mechanism is not the only one that produces a conflict-free system. Another example of conflict resolution mechanisms is the time-based conflict resolution mechanism used in the Mat\'ern type II model.

\sl The Mat\'ern type II construction \rm  gives each point $x$ a mark $t(x)$ which takes value in $[0,1]$ as an additional attribute of point $x$. This mark is interpreted as the time when the point `arrives' in the system. For convenience, we refer to this mark as the timer of the point. The spatial conflict between two ponts is resolved by a competition where whichever arrives first wins. Only the winning points belong to $\phi_{\rII}$, i.e.
\begin{align}
\phi_{\rII}=\;\{x \in \phi \st \textnormal{for all } y \in \phi, \bl c(x,y)=1 \,\Rightarrow\, t(x)<t(y) \}.
\end{align}
 In the literature this construction is sometimes referred to as the Mat\'ern hard-core model.
  
\sl The Mat\'ern type III construction \rm is proposed with the purpose of resolving conflicts while retaining as many points as possible. In this sense, it can be viewed as an improvement of the Mat\'ern type II model. The intuition behind this mechanism is as follows: when an object competes with others for space, it does not need to compete with those objects that have already been defeated. When the system $\phi$ contains only finitely many objects, we can give an explicit construction for the system $\phi_{\rIII}$. First, all the objects in $\phi$ are sorted in the increasing order of the values of their timers. Let $\{x_{i},\, i=1,2,\ldots\}$ be this ordering. We then construct an increasing sequence of sets $\{\phi_{\rIII}^{(i)},\, i=1,2,\ldots\}$ as follows. 
\begin{align*}
 &\phi_{\rIII}^{(1)}=\{x_{1}\};\\
&\phi_{\rIII}^{(i+1)}=\left\{\begin{array}{ll}
                             \phi_{\rIII}^{(i)} \cup \{x_{i+1}\} &\textnormal{ if }  c(x_{i+1},x_{j})=0 \bl \textnormal{for all } x_{j} \in \phi_{\rIII}^{(i)}, \\\noask
			     \phi_{\rIII}^{(i)} & \textnormal{ otherwise}.
                            \end{array} \right.
\end{align*}
The system $\phi_{\rIII}$ is defined as $\bigcup_{i=1}^{\infty}\phi_{\rIII}^{(i)}$. It is easily seen that $\phi_{\rIII}$ satisfies:
\begin{align}
 \phi_{\rm III}=\{x \in \phi \st \textnormal{for all } y \in \phi_{\rIII}, \bl c(x,y)=0 \Rightarrow t(x)<t(y)\}.
\end{align}
 In applications in physics, chemistry and material science, the Mat\'ern type III model is sometimes referred to as the Random Sequential Adsorption (RSA) model (see \cite{Pen02} and the citations herein).
 
 When $\phi$ contains infinitely many objects, there are configurations of $\phi$ such that the type III construction is not applicable. A simple example is the case where objects in $\phi$ are located at points in $\mathbb{Z}^{+}$ and $t(i)=i^{-1}$. We can easily see that there is no way to sort the objects in $\phi$ in the increasing order of their timers. Nevertheless, the construction of Mat\'ern type III model can still be extended to the case where $|\phi|=\infty$ under some mild condition. This construction is more involved and we postpone its discussion to Section \ref{S:type III extend}.

 Having recalled the classical constructions of the Mat\'ern models of type I, II and III, we can now introduce our extensions. As shown in Subsection \ref{SS: k to inf Matern}, these extensions form a bridge between Mat\'ern type II and Mat\'ern type III models.
\begin{enumerate}
 \item \sl The $0$-Mat\'ern and the $1$-Mat\'ern constructions \rm produce $\phi_{0}$, $\phi_{1}$ which are the original system $\phi$ and the system $\phi_{\rII}$ produced by the Mat\'ern type II model respectively.
 \item \sl The $2$-Mat\'ern construction \rm produces a subset $\phi_{2}$ of $\phi$ from $\phi_{1}$ in the following fashion. Each object $x$ in $\phi$ identifies in $\phi_{1}$ a set $A(x)$ of objects having a conflict with it. It is included in $\phi_{2}$ if it wins all the competitions with the elements of $A(x)$. In other words,  
\begin{align}
\nonumber A(x)&:=\;\{y \in \phi_{1} \st c(x,y)=1\};\\\noask
\nonumber \phi_{2}&:=\;\{x \in \phi \st \textnormal{for all } y \in A(x),\, t(x)<t(y) \}\\
		&=\;\{x \in \phi \st \textnormal{for all } y \in \phi_{1},\,c(x,y)=1 \Rightarrow t(x)<t(y) \}. 
\end{align}
 \item \sl The $(k+1)$-Mat\'ern construction \rm is recursively constructed from the $k$-Mat\'ern model in the same manner as the $2$-Mat\'ern model is constructed from the $1$-Mat\'ern model, i.e.\
 \begin{align}
\phi_{k+1}:=\;\{x \in \phi \st \textnormal{for all } y \in \phi_{k},\, c(x,y)=1 \Rightarrow t(x)<t(y) \}. 
\end{align}
\end{enumerate}

As a final note for this section, we want to stress that while Mat\'ern used timers taking values in $[0,1]$ (or more generally, in any bounded interval) in his construction, it is sometimes more convenient to have timers taking values in $\mathbb{R}^{+}$. This is in particular the case in Section \ref{S:Matern pgfl} where this extension allows us to obtain simpler analytical results. Hence, from now on we adopt this extension. 

\section{The Extended Mat\'ern Type III Construction}\label{S:type III extend}
Our aims here are to provide an extension of the Mat\'ern type III construction to the case where $\phi$ contains countably many points and to prove that the subset $\phi_{\infty}$ produced by this construction is a suitably defined \sl limit \rm of the sets $\phi_{k}$ when $k$ goes to infinity. We start by defining the \emph{conflict graph} associated with $\phi$ and $c$. Then we give our extension which is applicable only when the conflict graph has the \emph{finite history} property. The conflict graph and the finite history property are defined in Subsection \ref{SS: conflict graph}. We close this section by proving that when the later condition is satisfied, we have $\phi_{\infty}(\phi,t,c)=\bigcup_{k=0}^{\infty}\phi_{2k+1}=\bigcap_{k=0}^{\infty}\phi_{2k}$. For this reason, we refer to the construction in this section as the $\infty$-Mat\'ern model.

\subsection{Conflict Graph and $\infty$-Mat\'ern Model}\label{SS: conflict graph}
  For any two points $x$ and $y$ in $\phi$, we put a directed edge from $x$ to $y$ if $c(x,y)=1$ and $t(x)<t(y)$. Let $\mathcal{E}$ be the set of all such edges 
\[\mathcal{E}=\{(x,y) \in \phi^{2} \st c(x,y)=1 \textnormal{ and } t(x)<t(y)\}.\]
 The \emph{conflict graph} associated with $\phi$ and $c$ is the directed graph $\mathcal{G}=\{\phi,\,\mathcal{E}\}$. It is not difficult to check that $\mathcal{G}$ is an acyclic graph.

  We can now recursively define the order (i.e. an asymmetric, transitive binary relation)  $\leftarrowtail$ of $\phi$ as 
\begin{align}
 x \leftarrowtail y \textnormal{ if }\left\{\begin{array}{l}
							\textnormal{ either } (x,y) \in \mathcal{E},\\\noask
							\textnormal{ or } \textnormal{there exists } z \in \phi \st (z,y) \in \mathcal{E} \textnormal{ and } x\leftarrowtail z.
						      \end{array}\right.
\end{align}
 We call this the \emph{ancestor} order in $\phi$. For each $x$, let
\begin{align}
\mathcal{A}(x)=\{y \in \phi \st y \leftarrowtail x\} ,
\end{align}
be the set of its \sl ancestors \rm. $\mathcal{G}$ is said to have the \emph{finite history} property if $\mathcal{A}(x)$ is finite for all $x$ in $\phi$.

When $\mathcal{G}$ has the finite history property, let $e_{\infty}$ (the $\infty$ subscript is explained in Subsection \ref{SS: k to inf Matern}) be a $\{0,1\}$-value function taking elements of $\phi$ as argument and satisfies:
\begin{align}
 e_{\infty}(x)=\left\{\begin{array}{ll}
                          1 &\textnormal{ if } \mathcal{A}(x)=\emptyset,\\\noask
			  \displaystyle\prod_{y \in \phi,\, (y,x)\in \mathcal{E}}\,\bigl(1-e_{\infty}(y)\bigr) &\textnormal{ otherwise }.\label{E: infty Matern ind}
                         \end{array}\right.
\end{align}
We have then,
\begin{proposition}
 When the conflict graph $\mathcal{G}$ has the finite history property, there exists a unique function $e_{\infty}$ satisfying (\ref{E: infty Matern ind}).
\end{proposition}
\emph{Proof.} It is sufficient to show that $e_{\infty}(x)$ is uniquely determined for every $x$ in $\phi$. We do this by induction on $|\mathcal{A}(x)|$, this is possible by the finite history property assumption (so that $|\mathcal{A}(x)| <\infty$ a.s.). The base case is when $|\mathcal{A}(x)|=0$ so that $e_{\infty}(x)=1$ by definition. Note that such a $x$ always exists by the finite history  assumption (the argument is rather simple: if $|\mathcal{A}(x)|>0$ for all $x$ in $\phi$, we start with $x_{1}$ and build an infinite chain $\{x_{i},\,i=1,2,\ldots\}$ such that $x_{1}\leftarrowtail x_{2}\leftarrowtail x_{3}\cdots$\,. By the transitivity of $\leftarrowtail$, we have $\{x_{i},\,i=1,2,\ldots\} \subseteq \mathcal{A}(x_{1})$. Thus $|\mathcal{A}(x_{1})|=\infty$, which is a contradiction).

Now suppose that $e_{\infty}(x)$ is uniquely determined for every $x\in \phi$ such that $|\mathcal{A}(x)|<k$. Consider $x \in \phi$ such that $|\mathcal{A}(x)|=k$ (if such $x$ exists), then
\begin{align}
 e_{\infty}(x)=\prod_{y \in \phi,\, (y,x)\in \mathcal{E}}\bigl(1-e_{\infty}(y)\bigr),
\end{align}
by definition.
As $y \in \mathcal{A}(x)$ for all $y \in \phi$ such that $(y,x)\in \mathcal{E}$ and $|\mathcal{A}(x)|=k$ we have that the left-hand side is a product of finitely many terms. Moreover, by transitivity of $\leftarrowtail$, we know that $\mathcal{A}(y) \subset \mathcal{A}(x)$ for all $y \in \mathcal{A}(x)$. Hence, $|\mathcal{A}(y)|<|\mathcal{A}(x)|=k$ for all $y \in \mathcal{A}(x)$. Then, by the induction hypothesis, every $e_{\infty}(y)$ term appearing in the left-hand side is uniquely determined, so $e_{\infty}(x)$ is uniquely determined. \eop

The subset $\phi_{\infty}$ produced by the $\infty$-Mat\'ern model is defined as
\begin{align}
 \phi_{\infty}=\{x_{i} \in \phi \st e_{\infty}(x)=1\}.
\end{align}
It is easily checked that
\begin{align*}
 \phi_{\infty}=\{x_{i} \in \phi \st \textnormal{for all } y \in \phi_{\infty},
 \bl c(x,y)=0\, \Rightarrow \,t(x)<t(y)\},
\end{align*}
and that when $\phi$ contains finitely many points, $\phi_{\infty}$ and $\phi_{\rIII}$ are identical. This justifies our claim that the $\infty$-Mat\'ern model is the extension of Mat\'ern type III model to the case where $\phi$ may be countably infinite and its associated conflict graph has the finite history property.

\subsection{$\infty$-Mat\'ern 
Construction as the Limit of $k$-Mat\'ern Constructions} \label{SS: k to inf Matern} 
Let a sequence of functions $\{e_{k},\,k=1,2,\ldots\}$ which take elements of $\phi$ as argument be defined recursively by
\begin{align}
 e_{1}(x)&=\left\{\begin{array}{ll}
                      1 &\textnormal{ if } \mathcal{A}(x)=\emptyset\, ,\\\noask
		      0 &\textnormal{ otherwise}\,;
                     \end{array}\right.\label{E: 1-Matern ind}\\ 
 e_{k+1}(x)&=\prod_{y \in \phi,\, (y,x) \in \mathcal{E}}\bigl(1-e_{k}(y)\bigr), \label{E: k-Matern ind}
\end{align}
where the product taken over the empty set is $1$ by convention.

We first characterize the $k$-Mat\'ern model in terms of the $e_{k}$ function. 
\begin{proposition}
 For every $k \in \mathbb{Z}^{+}$
\begin{align*}
 \phi_{k}= \{x \in \phi \st e_{k}(x)=1\}.
\end{align*}
\end{proposition}
\emph{Proof}. We use induction on $k$. The base case $k=1$ is easily verified. Suppose that the result holds for some $k \geq 1$: $x \in \phi_{k}\, \Leftrightarrow \,e_{k}(x)=1$ for all $x$ in $\phi$. For $k+1$, we have
\begin{align*}
 \forall\, x \in \phi: x \in \phi_{k+1} &\Leftrightarrow \bigl(\forall\, y \in \phi_{k}: c(x,y)=1 \Rightarrow t(x)<t(y)\bigr) \\\noask
							  &\Leftrightarrow \bigl(\forall\, y\in \phi_{k}: (y,x)\notin \mathcal{E} \bigr)\; (\textnormal{by definition of } \mathcal{E})\\\noask
							  &\Leftrightarrow \bigl(\forall\, y\in \phi: y\in \phi_{k}\Rightarrow(y,x)\notin \mathcal{E} \bigr) \\\noask
							  &\Leftrightarrow \bigl(\forall\, y\in \phi: (y,x)\in \mathcal{E}\Rightarrow y\notin \phi_{k}\bigr) \\\noask
							  &\Leftrightarrow \bigl(\forall\, y\in \phi,(y,x)\in \mathcal{E}: y \notin \phi_{k}\bigr) \\\noask
							  &\Leftrightarrow \bigl(\forall\, y\in \phi,(y,x)\in \mathcal{E}: e_{k}(y)=0\bigr) \\\noask
							  &\Leftrightarrow \prod_{y\in \phi,\,(y,x)\in \mathcal{E}}\bigl(1-e_{k}(y)\bigr)=1\\\noask
							  &\Leftrightarrow e_{k+1}(x)=1\; (\textnormal{by definition})\,.
\end{align*}
\eop

In other words, the functions $\{e_{k},\,k=1,2,\ldots\}$ defined in (\ref{E: 1-Matern ind}) and (\ref{E: k-Matern ind}) are exactly the indicators that a point in $\phi$ belongs to the sets $\{\phi_{k},\,k=1,2,\ldots\}$ respectively . By simple induction arguments, we can also prove that
\begin{align*}
 &e_{2k+1}(x)\leq e_{2k}(x) \bl \textnormal{for all } k\in \mathbb{Z}^{+},x\in \phi,\\
 &e_{2k+1}(x)\leq e_{2k+3}(x) \bl \textnormal{for all } k\in \mathbb{Z}^{+},x\in \phi,\\
 &e_{2k+2}(x_{i})\leq e_{2k}(x) \bl \textnormal{for all } k\in \mathbb{Z}^{+},x\in \phi.
\end{align*}
So, we can deduce that $\phi_{1}\subseteq \phi_{3}\subseteq \cdots$ and $\phi_{0}\supseteq \phi_{2}\supseteq \cdots$, and hence
\begin{align}
 &\bigcup_{k=0}^{\infty} \phi_{2k+1} \subseteq \bigcap_{k=0}^{\infty} \phi_{2k}. \label{E: k-Matern monotonicity}
\end{align}
 Note that all the above observations hold even when the conflict graph does not have the finite history property. We now show that when the conflict graph does have the finite history property, $\bigcup_{k=0}^{\infty} \phi_{2k+1} = \bigcap_{k=0}^{\infty} \phi_{2k}=\phi_{\infty}$ and this justifies the $\infty$ subscript in $\phi_{\infty}$. To this end, the result of the next proposition is sufficient.
\begin{proposition}
 When the conflict graph $\mathcal{G}$ has the finite history property, for any point $x_{i}$ in $\phi$, 
\begin{align*}
 e_{\infty}(x)=\lim_{k\rightarrow \infty} e_{k}(x).
\end{align*}
\end{proposition}
\emph{Proof}. We prove by induction on $|\mathcal{A}(x)|$ that \[e_{|\mathcal{A}(x)|+1}(x)= e_{|\mathcal{A}(x)|+2}(x)=\cdots=e_{\infty}(x).\] The base case is when $|\mathcal{A}(x)|=0$. In this case $e_{1}(x)=e_{2}(x)=\cdots=e_{\infty}(x)=1$ by definition. Moreover
\begin{align*}
 e_{k}(x)= \prod_{x_{j}\in \phi, (y,x) \in \mathcal{E}}\bigl(1-e_{k-1}(y)\bigr)=1,
\end{align*}
   since $\{x_{j}\in \phi \st (x_{j},x_{i}) \in \mathcal{E}\} \subset \mathcal{A}(x_{i})=\emptyset$ and the product over an empty set is $1$ by convention.

Now suppose that the result holds for all $y$ such that$|\mathcal{A}(y)|\leq k$, for some $k \geq 0$ . Consider an $x$ such that $|\mathcal{A}(x)|=k+1$ (if such $x$ exists). We have for all $l$,
\begin{align*}
 e_{l}(x)=\prod_{y\in \phi, (y,x) \in \mathcal{E}}\bigl(1-e_{l-1}(y)\bigr).
\end{align*}
 Note that for all $y$ such that $(y,x) \in \mathcal{E}$, we must have $|\mathcal{A}(y)|\leq k$. By the induction hypothesis, $e_{l-1}(y)=e_{\infty}(y)$ for such $y$ and for $l>k+2$. Hence, for all $l>k+2$,
\begin{align*}
 e_{l}(x)&=\prod_{y\in \phi, (y,x) \in \mathcal{E}}\bigl(1-e_{l-1}(y)\bigr)=\prod_{y\in \phi, (y,x) \in \mathcal{E}}\bigl(1-e_{\infty}(y)\bigr)=e_{\infty}(x_{i}).
\end{align*}
 \eop

\section{Poisson Rain with Random Conflict Relation and its Mat\'ern Models}
\label{secprain}

\subsection{Poisson Rain with Random Conflict Relation}

\begin{definition}
A Poisson Rain with ground intensity $\Lambda$ is a PPP $\{(x,t(x))\}$ in $\mathbb{R}^{2}\times \mathbb{R}^{+}$ of intensity measure $\Lambda\times \mathcal{L}$ with $\mathcal{L}$ the Lebesgue measure.
\end{definition}
 Each 'point' in $\Phi$ is a pair $(x,t(x))$ with $x \in \mathbb{R}^{2}$ and $t(x) \in \mathbb{R}^{+}$. The $x$ component is understood as the position of a point and $t(x)$ is understood as the timer of the point. The name Poisson Rain stems from the interpretation of $\Phi$ as a collection of raindrops falling from the sky, the timer of a point is the time when it hits the ground and its position is the place where it does so. By abuse of notation, for each $x\in\mathbb{R}^{2}$, we write $x\in \Phi$ for ``there exists a $t$ such that $(x,t)\in \Phi$". So, when we refer to a point of $\Phi$ as a pair of location-timer, we use the pair notation $(x,t)$, and when we refer to it as a point in $\mathbb{R}^{2}$, we use the single element $x$ notation.
\begin{remark}
The Poisson Rain considered here is a special case of the extended Marked Point Process introduced in \cite[Definition 9.1.VI, p.~7]{DaVer08}. Note that the extended MPP in \cite{DaVer08} is used to construct the counting measure of a purely atomic random measure, where the mark of a point is used to represent the mass of the random measure at that point. As the mark of a point is used here to represent its arrival time, we use here the name Poisson Rain instead of extended marked Poisson Point Process.  
\end{remark} 
Let $h$ be a function satisfying \eqref{E: finite conflict}.
The random conflict relation is extended to the Poisson Rain as follows:
 
\begin{definition}
A Poisson Rain with Random Conflict Relation (PRRCR) with ground intensity $\Lambda$ and expected conflict function $h$ is a pair $(\Phi,C)$ where $\Phi$ is a Poisson Rain with ground intensity $\Lambda$ and $C=\{C(x,y),\,x,y \in \Phi\}$ is a family of $\{0,1\}$ value r.v.s indexed by unordered pairs of locations in $\Phi$ satisfying
\begin{itemize}
\item[1] $C$ is non-reflexive and symmetrical a.s.; 
\item[2] given a realization of $\Phi$, $C$ is a family of independent r.v.s with the exceptions given by condition 1; and
\item[3] $\mathbb{P}\bigl(C(x,y)\mid \Phi \mbox{ and } x,y \in \Phi\bigr)=h(x,y)$.
\end{itemize}  
\end{definition}

For completeness, we provide here a construction of a PRRCR. Let $\{\Psi_{i},\, i=1,2,\ldots\}$ be a family of i.i.d.\ MPPPs with i.i.d.\ marks of ground intensity $\Lambda$. Each point $x$ in $\Psi_{i}$ is equipped with a mark $\textbf{u}(x)=\{\tau(x),u_{j,k}(x),\,j,k=1,2,\ldots\}$ which is a family of i.i.d.\ r.v.s uniformly distributed in $[0,1]$. Let 
\begin{align*}
\Psi^{'}_{i}:=\{\bigl(x,\tau(x)+i,\{u_{j,k}(x)\}\bigr) \mbox{ for all } \bigl(x,\tau(x),\{u_{j,k}(x)\}\bigr) \in \Psi_{i}\}.
\end{align*} The Poisson Rain $\Phi$ is defined as $\bigcup_{i=1}^{\infty}\Psi^{'}_{i}$. To determine the random conflict relation, for each $x$ in $\Phi$, we number the points in $\Psi'_{i}$ in the increasing order of their distance to the centre of the space and associate to the $j^{th}$ point in this numbering the r.v.\ $u_{i,j}(x)$ in the mark of $x$. Now, consider any two points $x$, $y$ in $\Phi$. Let $u_{i,j}(x)$ be the r.v.\ corresponding to $y$ in the mark of $x$ and $u_{k,l}(y)$ be the r.v.\ corresponding to $x$ in the mark of $y$, we let $C(x,y)=1$ iff $\min\bigl(u_{k,l}(y),u_{k,l}(y)\bigr)>1-\sqrt{h(x,y)}$. It is then easily verified that $(\Phi,C)$ defined in this manner is indeed a PRRCR with ground intensity $\Lambda$ and expected conflict function $h$.

\subsection{$\infty$-Mat\'ern Construction of Poisson Rain with Random Conflict Relations} 
Throughout this paper, we assume that the measure $\Lambda$ and the function $h$ satisfy
\begin{align}
 \sup_{x \in \mathbb{R}^{d}} \int_{\mathbb{R}^{d}}h(x,y)dy = \mathcal{N}<\infty. \label{E: finite conflict}
\end{align}
 The main result of this section is
\begin{proposition}\label{P: well def cond}
 For any PRRCR $\Phi$ of ground intensity measure $\Lambda$ and expected conflict function $h$ satisfying (\ref{E: finite conflict}), its corresponding conflict graph has the finite history property ( see definition in Subsection \ref{SS: conflict graph}) a.s.   
\end{proposition}
\emph{Proof.} Consider a typical point $x$ in $\Phi$. Recall that $\mathcal{G}=(\Phi,\mathcal{E})$ is the conflict graph and $\mathcal{A}(x)$ is the set of ancestors of $x$ in the conflict graph. It is sufficient to prove that $\mathcal{A}(x)$ is finite for every $x$. Let 
\begin{align*}
\mathcal{A}^{(l)}(x)=& \{y \in \Phi \st \textnormal{exists } x_{0},x_{1},\ldots,x_{l} \in \Phi \\
&\st x_{0}=x, x_{l}=y \textnormal{ and }  (x_{1},x_{0}),\ldots,(x_{l},x_{1-1}) \in \mathcal{E}\}.
\end{align*}
$\mathcal{A}(x)$ can be rewritten as:
\begin{align*}
 \mathcal{A}(x)=\bigcup_{l=1}^{\infty}\mathcal{A}^{(l)}(x).
\end{align*}
We now prove  that $\textbf{E}\left[|\mathcal{A}^{(l)}(x)\right]\leq \frac{t(x)^{l}\mathcal{N}^{l}}{l!}$. First notice that
\begin{align*}
 |\mathcal{A}^{(l)}(x)|&\leq\, \sum_{x_{1},\ldots,x_{l} \in \Phi}\textbf{1}_{(x_{1},x)\in \mathcal{E}}\prod_{k=2}^{l}\textbf{1}_{(x_{k},x_{k-1})\in \mathcal{E}}\\\noask
			   &=\, \sum_{x_{1},\ldots,x_{l} \in \Phi,\,x_{0}=x}\prod_{k=0}^{l-1}C(x_{k},x_{k+1})\textbf{1}_{t(x_{k})>t(x_{k+1})}  .
\end{align*}
We consider $\Phi$ as a PPP of intensity $\Lambda\times\mathcal{L}$ in $\mathbb{R}^{d}\times \mathbb{R}^{+}$, with $\mathcal{L}$ Lebesgue measure, and apply the multivariate Campbell formula to get
\begin{align*}
 \mathbb{E}\left[|\mathcal{A}^{(l)}(x)|\right] &\leq& \mathbb{E}\left[\, \sum_{x_{0},\ldots,x_{l} \in \Phi,\,x_{0}=x}\prod_{k=0}^{l-1}C(x_{k},x_{k+1})\textbf{1}_{t(x_{k})>T(x_{k+1})} \right]\\\noask
 &=&\int_{\left(\mathbb{R}^{d}\right)^{l}}\int_{\left(\mathbb{R}^{+}\right)^{l}}\;\textbf{1}_{t(x)>t_{1}>t_{2}>\cdots>t_{l}}\;h(x,x_{1})\;\prod_{k=1}^{l-1}h(x_{l},x_{l+1})\;\\\noask
& &dt_{1}\;\ldots\; dt_{l}\;\Lambda(dx_{1})\;\ldots\;\Lambda(dx_{l})\\\noask
&=&\int_{\left(\mathbb{R}^{d}\right)^{l}}\;\frac{t(x)^{l}}{l!}\;h(x,x_{1})\;\prod_{k=1}^{l-1}h(x_{k},x_{k+1})\;\Lambda(dx_{1})\;\ldots\;\Lambda(dx_{l}).
\end{align*}
 As $h$ and $\Lambda$ satisfy the condition (\ref{E: finite conflict}), we deduce that
\begin{align*}
 \mathbb{E}\left[|\mathcal{A}^{(l)}(x)|\right] &\leq\; \int_{\left(\mathbb{R}^{d}\right)^{l}}\;\frac{t(x)^{l}}{l!}\;h(x,x_{1})\;\prod_{k=1}^{l-1}h(x_{l},x_{l+1})\;\Lambda(dx_{1})\;\ldots\;\Lambda(dx_{l})\\
						   &\leq\; \frac{t(x)^{l}}{l!}\mathcal{N}^{l}.
\end{align*}
So
\begin{align*}
 \mathbb{E}\left[|\mathcal{A}(x)|\right] \leq \sum_{l=1}^{\infty} \mathbb{E}\left[|\mathcal{A}^{(l)}(x)|\right] \leq\; \sum_{l=1}^{\infty} \frac{t(x)^{l}}{l!}\mathcal{N}^{l}= \exp\{t(x) \mathcal{N}\}<\infty.
\end{align*}
This implies that $|\mathcal{A}(x)|<\infty$ a.s. 
\eop

In particular, Proposition \ref{P: well def cond} shows that when a PRRCR satisfies \eqref{E: finite conflict}, its $\infty$-Mat\'ern model is well defined and is the \emph{limit} of the $k$-Mat\'ern models.
\subsection{Timer-based Restrictions and $k$-Mat\'ern Models}
	As we will see in the next sections, to study the distribution of the PPs induced by the Mat\'ern models of the PRRCRs, it is useful to divide these PPs to thin layers such that each layer looks like a Poisson PP. Such a partition is defined in this subsection. 
    Given a PRRCR $\Phi$, we define the restriction $T_{s,t}$ to the interval $[s,t)$ as 
      \begin{align}
      T_{s,t}(\Phi)= \{x \in \Phi \st t(x) \in [s,t)\}.
      \end{align}
      This restricted version of $\Phi$ inherits the natural conflict relation  from $\Phi$. When $s=0$, the above notation is reduced to $T_{t}(\Phi)$. Consider $T_{t}(\Phi)$ as a PRRCR,we can easily verify that the $k$-Mat\'ern  construction is applicable to $T_{s,t}(\Phi)$ and the $\infty$-Mat\'ern type construction is also applicable a.s. given that $\Lambda$ and $h$ satisfy \eqref{E: finite conflict}. In particular, let the constructed PPs inherit the natural conflict relation from their original PRRCR. The following facts can be easily proved:
\begin{itemize}
     \item $T_{t}(\mathcal{M}_{j}(\Phi))=\mathcal{M}_{j}(T_{t}(\Phi))$ for $j= 0,1,\ldots\infty$;
     \item $\mathcal{M}_{j}(\Phi)=\{x \st C(x,y)= 0\; \forall \;y \in T_{t(x)}(\mathcal{M}_{j-1}\Phi)\}$ for $j= 0,1,\ldots$; and
     \item $\mathcal{M}_{\infty}(\Phi)=\{x \st C(x,y)= 0\; \forall \;y \in T_{t(x)}(\mathcal{M}_{\infty}(\Phi))\}$.
\end{itemize} 

In other words, the first claim asserts that the restriction to $[0,t)$ of the Mat\'ern models of $\Phi$ are the same as the Mat\'ern models of the restriction to $[0,t)$ of $\Phi$ while the two other claims are just reformulations of the Mat\'ern models definitions. These three claims are based on the fact that for $j=1,2,\ldots,\infty$, the event $x \in \mathcal{M}_{j}(\Phi)$  depends only on the realization of the points in $\Phi$ whose timers are smaller than $t(x)$.

We conclude with a result that will be used frequently in the next section. It allows us to approximate the PPs $T_{s,t}\bigl(\mathcal{M}_{j}(\Phi)\bigr), j=1,2,\ldots,\infty$ by other PPs which have much simpler structures.

\begin{proposition}\label{P:inf thin layer}
 For every realization of the PRRCR $(\Phi,C)$ such that its associated conflict graph has the finite history property and for every $0<s<t$,
\begin{align}
\nonumber&\Delta_{j,d,s,t} \subseteq T_{s,t}\bigl(\mathcal{M}_{j}(\Phi)\bigr)  \subseteq \Delta_{j,u,s,t} \textnormal{ for $j=1,2,\ldots$};\\
&\Delta_{\infty,d,s,t} \subseteq T_{s,t}\bigl(\mathcal{M}_{\infty}(\Phi)\bigr) \subseteq \Delta_{\infty,u,s,t},
\end{align}
where
\begin{align}
\nonumber& \Delta_{j,d,s,t}=\{x \in T_{s,t}(\Phi) \st C(x,y)=0 \;\forall\; y \in T_{s}\bigl(\mathcal{M}_{j-1}(\Phi)\bigr)\cup T_{s,t}(\Phi)\};\\
\nonumber& \Delta_{j,u,s,t}=\{x \in T_{s,t}(\Phi) \st C(x,y)=0 \;\forall\; y \in T_{s}\bigl(\mathcal{M}_{j-1}(\Phi)\bigr)\};\\ 
\nonumber& \Delta_{\infty,d,s,t}=\{x \in T_{s,t}(\Phi) \st C(x,y)=0 \;\forall\; y \in  T_{s}\bigl(\mathcal{M}_{\infty}(\Phi)\bigr) \cup T_{s,t}(\Phi)\};\\ 
& \Delta_{\infty,u,s,t}=\{x \in T_{s,t}(\Phi) \st C(x,y)=0 \;\forall\; y \in  T_{s}\bigl(\mathcal{M}_{\infty}(\Phi)\bigr)\}. 
\end{align}
\end{proposition}
\emph{Proof.}
 Since $\mathcal{M}_{j}(\Phi)=\{x \st C(x,y)= 0\; \forall \;y \in T_{t(x)}\bigl(\mathcal{M}_{j-1}(\Phi)\bigr)\}$
and $T_{s}(\Phi)\subseteq T_{t(x)}(\Phi)\subseteq T_{t}(\Phi)$ for every $x \in T_{s,t}(\Phi)$, the first double inclusion follows directly. The second one is proved similarly using the third claim above and taking into account that $T_{s}\bigl(\mathcal{M}_{\infty}(\Phi)\bigr)\subseteq T_{t(x)}\bigl(\mathcal{M}_{\infty}(\Phi)\bigr)\subseteq T_{s}\bigl(\mathcal{M}_{\infty}(\Phi)\bigr)\cup T_{s,t}(\Phi)$ for every $x\in T_{s,t}(\Phi)$.
\eop

\section{Probability Generating Functionals}\label{S:Matern pgfl}
This section gathers the main results of this paper regarding the evolution of the p.g.fls of the $k$-Mat\'ern models and the $\infty$-Mat\'ern model . Even though it is harder to define, the p.g.fls of the $\infty$-Mat\'ern model is easier to work with. For this reason, we give the results concerning these p.g.fls in Subsection \ref{SS: inf Matern} after some background on p.g.fls of PPs in Subsection \ref{SS:p.g.fl}. Then, the method is extended to study the p.g.fls of the $k$-Mat\'ern models in Subsection \ref{SS: k Matern}.
\subsection{Preliminaries}\label{SS:p.g.fl}
For any PP $\Xi$, we define for each function $v$ taking value in $[0,1]$
\begin{align}
G_{\Xi}(v)=\mathbb{E}\left [\prod_{x \in \Xi}v(x)\right ].
\end{align}
This is called the probability generating functional (p.g.fl) of $\Xi$ at $v$. For the p.g.fl to be well-defined and non-trivial, we need that 
$\left |\sum_{x \in \Xi}\log\bigl(v(x)\bigr)\right |< \infty$ a.s. As $\left |\log\bigl(v(x)\bigr)\right |\leq \bigl(1-v(x)\bigr)$, we consider only the p.g.fl of $\Xi$ at functions $v$ satisfying 
\begin{align}
\mathbb{E}\left [\sum_{x \in \Xi}\bigl(1-v(x)\bigr)\right ]=\int_{\mathbb{R}^{d}}\bigl(1-v(x)\bigr)m_{\Xi}(dx) < \infty, \label{E: gen def cond}
\end{align}
with $m_{\Xi}(\cdot)$ the intensity (first moment) measure of $\Xi$. 
An important special case is the p.g.fls of PPPs, which can be computed in closed forms \cite{DaVer08},\cite{BB09}.
    \begin{theorem}\label{T:PoissLaplace}
    Let $\Psi$ be a PPP of intensity measure $\Lambda$ and $v$ is a function taking value in $[0,1]$ such that $\int_{\mathbb{R}^{d}}\bigl(1-v(x)\bigr)\,\Lambda(dx)< \infty$. Then,
    \begin{align}
    G_{\Psi}(v)=\exp\left \{-\int_{\mathbb{R}^{d}}\bigl(1-v(x)\bigr)\,\Lambda(dx)\right \}.
    \end{align}
    \end{theorem} 

\subsection{$\infty$-Mat\'ern Models}\label{SS: inf Matern}
We are interested in the functionals $f_{k,\Lambda}(t,v):=G_{T_{0,t}(\mathcal{M}_{k}(\Phi))}(v(\cdot))$ with $t \geq 0$, where $\Phi$ is a PRRCR with ground intensity $\Lambda$, expected conflict function $h$  and $v$ is a function from $\mathbb{R}^{2}$ to $[0,1]$ such that:
\begin{align}
 \int_{\mathbb{R}^{2}}\bigl(1-v(x)\bigr)\,\Lambda(dx) <\infty. \label{E: Pois gen func cond}
\end{align}
In the other words, for each fixed $t$ this gives us the p.g.fls of the PPs induced by Mat\'ern models restricted to $[0,t]$.
\paragraph{Remarks}
\begin{itemize}
 \item The condition (\ref{E: Pois gen func cond}) is equivalent to
$\mathbb{E}\left[\sum_{x \in T_{0,t}(\Phi)}\bigl(1-v(x)\bigr)\right]<\infty$ (Campbell's formula) for every $t\geq 0$. Thus, this condition guarantees
that the p.g.fl of $\Phi$ is well-defined and non-trivial at $v$. In particular, since 
\[T_{0,t}(\mathcal{M}_{k}(\Phi))\subset T_{0,t}(\Phi) \mbox{ a.s. },\]
we have $\mathbb{E}\left[\sum_{x \in T_{0,t}(\mathcal{M}_{k}(\Phi))}\bigl(1-v(x)\bigr)\right]<\infty$, so that $G_{T_{0,t}(\mathcal{M}_{k}(\Phi))}$ is well-defined at $v$.
\item Although there are functions $v$ such that $G_{T_{0,t}(\mathcal{M}_{k}(\Phi))}$
is well-defined at $v$ while (\ref{E: Pois gen func cond}) is not satisfied, we do not consider 
this case here for two reasons:
\begin{itemize}
\item All information about the distribution of a PP can be extracted from its p.g.fl at functions $v$ having bounded support. These functions satisfy    
(\ref{E: Pois gen func cond}) automatically.
\item As we see in Subsections \ref{SS: inf Matern} and \ref{SS: k Matern}, this condition guarantees nice convergence properties.
\end{itemize}
\end{itemize}
If not otherwise stated, we assume that every function considered here satisfies \eqref{E: Pois gen func cond}.
Fix such a function $v$, we first show that $f_{\infty,\Lambda}(t,v)$ is continuous in $t$. 
\begin{proposition}\label{P:f cont}
For every function $v$ satisfying condition (\ref{E: Pois gen func cond}),
$f_{\infty,\Lambda}(t,v)$ is continuous in $t$.
\end{proposition}
\emph{Proof.}
For every $t$ and $\epsilon$ positive
\begin{align*}
 f_{\infty,\Lambda}(t+\epsilon,v)&=\mathbb{E}\left[\prod_{x \in T_{t+\epsilon}(\mathcal{M}_{\infty}(\Phi))}v(x)\right]\\\noask 
&=\mathbb{E}\left[\prod_{x \in T_{t}(\mathcal{M}_{\infty}(\Phi))}v(x)\prod_{x \in T_{t,t+\epsilon}(\mathcal{M}_{\infty}(\Phi))}v(x)\right].
\end{align*}
  As $T_{t,t+\epsilon}\bigl(\mathcal{M}_{\infty}(\Phi)\bigr) \subseteq T_{t,t+\epsilon}(\Phi)$ a.s., we deduce that
 \begin{align*}
&1\geq \prod_{x \in T_{t,t+\epsilon}(\mathcal{M}_{\infty}(\Phi))}v(x) \geq \prod_{x \in T_{t,t+\epsilon}(\Phi)}v(x) \mbox{ a.s.}
\end{align*}
As $T_{t,t+\epsilon}$ is independent with $T_{t}\bigl(\mathcal{M}_{\infty}(\Phi)\bigr)$,
 \begin{align*}
f_{\infty,\Lambda}(t,v)e^{-\epsilon\int_{\mathbb{R}^{2}}(1-v(x))\Lambda(dx)}
&=\; \mathbb{E}\left[\prod_{x \in T_{t}(\mathcal{M}_{\infty}(\Phi))}v(x)\prod_{x \in T_{t,t+\epsilon}(\Phi)}v(x)\right]\\\noask
&\leq\; \mathbb{E}\left[\prod_{x \in  T_{t+\epsilon}(\mathcal{M}_{\infty}(\Phi))}v(x)\right]\\\noask
&=\; f_{\infty,\Lambda}(t+\epsilon,v) \leq f_{\infty,\Lambda}(t,v),
\end{align*}
where the inequality comes from the fact that $v\leq 1$.
Following the same method, we get
\begin{align*}
&f_{\infty,\Lambda}(t-\epsilon,v)e^{-\epsilon\int_{\mathbb{R}^{2}}|1-v(x)|\,\Lambda(dx)}\leq f_{\infty,\Lambda}(t,v)\leq f_{\infty,\Lambda}(t-\epsilon,v).
\end{align*}
Letting $\epsilon$ go to $0$ completes this proof.
\eop

Let $H$ be the mapping that associates to a function $v: {\mathbb{R}^{2}}\to [0,1]$ and a point
$x\in \mathbb{R}^{2}$ the function 
\begin{align}
H(v,x):\,y \longmapsto v(y)\bigl(1-h(x,y)\bigr) \label{E:h transform}
\end{align}
from ${\mathbb{R}^{2}}$ to $[0,1]$.
\begin{theorem}\label{T:CSMA eq}
For any locally finite measure $\Lambda$, the functional $f_{\infty,\Lambda}$ satisfies the following system of equations,
\begin{eqnarray}
\nonumber f_{\infty,\Lambda}(0,v)&=&1;\\\noask
 \frac{df_{\infty,\Lambda}(t,v)}{dt}&=& -\int_{\mathbb{R}^{2}}f_{\infty,\Lambda}\bigl(t,H(v,x)\bigr)\bigl(1-v(x)\bigr)\Lambda(dx).\label{E:CSMA eq}
\end{eqnarray}
\end{theorem}
The main idea behind this theorem is to divide $\mathcal{M}_{\infty}(\Phi)$ into very thin layers $T_{t,t+\epsilon}\bigl(\mathcal{M}_{\infty}(\Phi)\bigr)$. For $\epsilon>0$ small enough, the points in each layer are so sparse that there is almost no conflict between them. Then we can consider each layer as a PPP. \\
In particular, we need to prove that
\begin{align}
&\lim_{\epsilon\rightarrow \infty}\frac{f_{\infty,\Lambda}(t+\epsilon,v)-f_{\infty,\Lambda}(t,v)}{\epsilon}
=-\int_{\mathbb{R}^{2}}f_{\infty,\Lambda}\bigl(t,H(v,x)\bigr)\bigl(1-v(x)\bigr)\Lambda(dx); \label{E:right limit}\\\noask
&\lim_{\epsilon\rightarrow \infty}\frac{f_{\infty,\Lambda}(t,v)-f_{\infty,\Lambda}(t-\epsilon,v)}{\epsilon} 
=-\int_{\mathbb{R}^{2}}f_{\infty,\Lambda}\bigl(t,H(v,x)\bigr)\bigl(1-v(x)\bigr)\Lambda(dx)\label{E:left limit}.
\end{align}

For this purpose, we need the following lemmas.
\begin{lemma}\label{L: main eq}
 Let $\Xi$ be a PP and $v$ is a function taking value in $[0,1]$ such that
\begin{align*}
 \mathbb{E}\left[\sum_{x \in \Xi}\bigl(1-v(x)\bigr)\right]<\infty.
\end{align*}
We have
\begin{align}
 \prod_{x  \in \Xi}v(x)= 1+\sum_{i=1}^{\infty}(-1)^{i}\sum_{(x_{1},\cdots,x_{i}) \in \Xi^{i!}}\prod_{j=1}^{i}\bigl(1-v(x_{j})\bigr) \textnormal{ a.s.},\label{E: main eq}
\end{align}
where $\Xi^{(i!)}$ is the set of unordered $i$-tuples of mutually different points in $\Xi$.
\end{lemma}
\emph{Proof.}
Since $\sum_{x \in \Xi}\log\bigl(v(x)\bigr) \leq \sum_{x \in \Xi}\bigl(1-v(x)\bigr)<\infty$ a.s., 
\[\prod_{x  \in \Xi}v(x)= \exp\left \{\sum_{x \in \Xi}\log\bigl(v(x)\bigr)\right \}\]
 is well-defined and is finite a.s.
 
 Now we need to prove that the series in the right hand side of (\ref{E: main eq}) converges. For this purpose, it is sufficient to show that
\begin{align*}
  1+\sum_{i=1}^{\infty}\sum_{(x_{1},\cdots,x_{i}) \in \Xi^{i!}}\prod_{j=1}^{i}\bigl(1-v(x_{j})\bigr)<\infty \textnormal{ a.s.}
\end{align*}
 Note that
\begin{align*}
 \sum_{(x_{1},\cdots,x_{i}) \in \Xi^{i!}}\prod_{j=1}^{i}\bigl(1-v(x_{j})\bigr)&=\frac{1}{i!}\sum_{x_{1},\cdots,x_{i} \textnormal{ mutually different} \in \Xi}\textnormal{ }\prod_{j=1}^{i}\bigl(1-v(x_{j})\bigr)\\
&\leq \frac{1}{i!}\left(\sum_{x \in \Xi}\bigl(1-v(x)\bigr)\right)^{i}.
\end{align*}
Hence,
\begin{align*}
  1+\sum_{i=1}^{\infty}\sum_{(x_{1},\cdots,x_{i}) \in \Xi^{i!}}\prod_{j=1}^{i}\bigl(1-v(x_{j})\bigr)&\leq 1+\sum_{i=1}^{\infty}\frac{1}{i!}\left(\sum_{x \in \Xi}\bigl(1-v(x)\bigr)\right)^{i}\\
&=\exp\left \{\left(\sum_{x \in \Xi}\bigl(1-v(x)\bigr)\right)\right \}<\infty.
\end{align*}
The equality can now be obtained by writing
\begin{align*}
 \prod_{x\in \Xi}v(x)= \prod_{x \in \Xi}\Bigl(1-\bigl(1-v(x)\bigr)\Bigr).
\end{align*}
\eop

\begin{lemma}\label{L:limit 0}
Let $(\Phi,C)$ be a PRRCR with ground intensity $\Lambda$ and expected conflict function $h$. For every $t\geq 0$ and $\epsilon>0$ small enough,
\begin{align}
\nonumber&\left|\mathbb{E}\left[\prod_{x \in T_{t,t+\epsilon}(\mathcal{M}_{\infty}(\Phi))}v(x) \,\Bigg|\, \Phi_{t}\right]-1+\epsilon\int_{\mathbb{R}^{2}}\bigl(1-v(x)\bigr)\,\Lambda_{\infty,t}(dx)\right|\\
&\leq \epsilon^{2}\left(2\left(\int_{\mathbb{R}^{2}}\bigl(1-v(x)\bigr)\,\Lambda(dx)\right)^{2}+\overline{H}\int_{\mathbb{R}^{2}}\bigl(1-v(x)\bigr)\,\Lambda(dx)\right) \mbox{ a.s. },
\end{align}
where $\overline{H}$ is defined in (\ref{E: finite conflict}) and $\Lambda_{\infty,t}$ is the random measure in $\mathbb{R}^{2}$ satisfying
\begin{align}
 \Lambda_{\infty,t}(dx)=\prod_{y \in T_{t}(\mathcal{M}_{\infty}(\Phi))}\bigl(1-h(x,y)\bigr)\,\Lambda(dx).
\end{align}
\end{lemma}
\emph{Proof.} See Appendix \ref{A: Proof Lemma Limit 0}.

 Note that given $\Phi_{t}$, the p.g.fl of $\Delta_{\infty,u,t,t+\epsilon}$, which is a PPP of intensity measure $\epsilon\Lambda_{\infty,t}$, is
\begin{align*}
 G_{\Delta_{\infty,u,t,t+\epsilon}}(v)=\exp\left\{-\epsilon\int_{\mathbb{R}^{2}}\bigl(1-v(x)\bigr)\Lambda_{\infty,t}(dx)\right\}.
\end{align*}
Hence,
\begin{align*}
 \lim_{\epsilon \rightarrow 0}\frac{G_{\Delta_{\infty,u,t,t+\epsilon}}(v)-1}{\epsilon}&=-\int_{\mathbb{R}^{2}}\bigl(1-v(x)\bigr)\Lambda_{\infty,t}(dx)\\\noask
&= \lim_{\epsilon \rightarrow 0}\frac{G_{T_{t,t+\epsilon}(\mathcal{M}_{\infty}(\Phi))}(v)-1}{\epsilon}.
\end{align*}
Thus, Lemma \ref{L:limit 0} justifies our intuition that when the time scale is small, the effect of conflict is negligible, and we can regard the thin layer $T_{t,t+\epsilon}\bigl(\mathcal{M}_{\infty}(\Phi)\bigr)$ as a PPP. Such property, which we call the \emph{quasi-Poisson} property, plays an important role in the subsequent studies.

Now we can proceed to the proof of Theorem \ref{T:CSMA eq}, where  $\infty$ and $\Lambda$ subscripts are dropped to avoid cumbersome notation. 
Note that for $r>s$
\begin{eqnarray*}
 f(r,v)-f(s,v) = \mathbb{E}\left[\prod_{x \in T_{s}(\mathcal{M}_{\infty}(\Phi))}v(x)
\left(\prod_{y \in T_{s,r}(\mathcal{M}_{\infty}(\Phi))}v(y)-1\right)\right].
\end{eqnarray*}
In order to evaluate the last expression, we need the conditional probability
\begin{align}
 \mathbb{E}\left[\left.\prod_{y \in T_{s,r}(\mathcal{M}_{\infty}(\Phi))}v(y)-1\,\right|\,\Phi_{s}\right]. \label{E: Proof eq 2}
\end{align}
Put $s=t$ and $r=t+\epsilon$, by Lemma \ref{L:limit 0},  
\begin{align*}
&\left|\left(\mathbb{E}\left[\left .\prod_{y \in T_{t,t+\epsilon}(\mathcal{M}_{\infty}(\Phi))}v(y)\,\right |\,\Phi_{t}\right]-1\right)\epsilon^{-1}+\int_{\mathbb{R}^{2}}\bigl(1-v(x)\bigr)\,\Lambda_{\infty,t}(dx)\right|\\\noask
&\leq \epsilon\left(2\left(\int_{\mathbb{R}^{2}}\bigl(1-v(x)\bigr)\,\Lambda(dx)\right)^{2}+\overline{H}\,\int_{\mathbb{R}^{2}}\bigl(1-v(x)\bigr)\,\Lambda(dx)\right) \mbox{ a.s. }
\end{align*}
Let $\mathfrak{C}=2\left(\int_{\mathbb{R}^{2}}\bigl(1-v(x)\bigr)\,\Lambda(dx)\right)^{2}+\overline{H}\,\int_{\mathbb{R}^{2}}\bigl(1-v(x)\bigr)\,\Lambda(dx)$, we have
\begin{align}
\nonumber&\left|\frac{f(t+\epsilon,v)-f(t,v)}{\epsilon}+\mathbb{E}\left[\prod_{y \in T_{t}(\mathcal{M}_{\infty}(\Phi))}v(y)\left(\int_{\mathbb{R}^{2}}\bigl(1-v(x)\bigr)\,\Lambda_{\infty,t}(dx)\right)\right]\right|\\\noask
&\leq \epsilon \mathfrak{C}\mathbb{E}\left[\prod_{y \in T_{t}(\mathcal{M}_{\infty}(\Phi))}v(y)\right]=\epsilon \mathfrak{C} f(t,v)\leq \epsilon\mathfrak{C}.\label{E: Proof eq 1}
\end{align}
 Moreover,
\begin{align*}
 &\mathbb{E}\left[\prod_{y \in T_{t}(\mathcal{M}_{\infty}(\Phi))}v(y)\left(\int_{\mathbb{R}^{2}}\bigl(1-v(x)\bigr)\,\Lambda_{\infty,t}(dx)\right)\right]\\\noask
 &=\mathbb{E}\left[\prod_{y \in T_{t}(\mathcal{M}_{\infty}(\Phi))}v(y)\left (\int_{\mathbb{R}^{2}}\bigl(1-v(x)\bigr)\,\left (\prod_{y \in T_{t}(\mathcal{M}_{\infty}(\Phi))}\bigl(1-h(x,y)\bigr)\right )\Lambda(dx)\right )\right]\\\noask
 &=\mathbb{E}\left[\int_{\mathbb{R}^{2}}\bigl(1-v(x)\bigr)\left (\prod_{y \in T_{t}(\mathcal{M}_{\infty}(\Phi))}v(y)\bigl(1-h(x,y)\bigr)\,\right )\Lambda(dx)\right]\\\noask
  &=\mathbb{E}\left[\int_{\mathbb{R}^{2}}\bigl(1-v(x)\bigr)\left (\prod_{y \in T_{t}(\mathcal{M}_{\infty}(\Phi))}H(v,x)(y)\right )\,\Lambda(dx)\right].
\end{align*}

As the term inside of the integration is a positive r.v.\, we can change the order of expectation and integration,
\begin{align*}
&\mathbb{E}\left[\int_{\mathbb{R}^{2}}\bigl(1-v(x)\bigr)\prod_{y \in T_{t}(\mathcal{M}_{\infty}(\Phi))}H(v,x)(y)\,\Lambda(dx)\right]\\\noask
&=\int_{\mathbb{R}^{2}}\bigl(1-v(x)\bigr)f\bigl(t,H(v,x)\bigr)\,\Lambda(dx).
\end{align*}

Letting $\epsilon$ goes to $0$ in (\ref{E: Proof eq 1}) gives us (\ref{E:right limit}). To prove (\ref{E:left limit}), we put $s=t-\epsilon$ and $r=t$ in (\ref{E: Proof eq 2}). Proceeding as above, we obtain
\begin{align*} 
&\left|\frac{f(t,v)-f(t-\epsilon,v)}{\epsilon}+\mathbb{E}\left[\prod_{y \in T_{t-\epsilon}(\mathcal{M}_{\infty}(\Phi))}v(y)\left(\int_{\mathbb{R}^{2}}\bigl(1-v(x)\bigr)\Lambda_{\infty,t-\epsilon}(dx)\right)\right]\right|\\\noask
&\leq \epsilon \mathfrak{C} f(t-\epsilon,v)\leq \epsilon \mathfrak{C}.
\end{align*}
We can also show that
\begin{align*}
&\mathbb{E}\left[\prod_{y \in T_{t-\epsilon}(\mathcal{M}_{\infty}(\Phi))}v(y)\left(\int_{\mathbb{R}^{2}}\bigl(1-v(x)\bigr)\Lambda_{\infty,t-\epsilon}(dx)\right)\right]\\\noask
&=\int_{\mathbb{R}^{2}}(1-v(x))f(t-\epsilon,H(v,x)))\Lambda(dx).
\end{align*}
We then get (\ref{E:left limit}) by letting $\epsilon$ goes to $0$ and using the continuity of $f$ in $t$ (Proposition \ref{P:f cont}).
\eop

We can also obtain an upper bound and a lower bound on $f_{\infty,\Lambda}$:
\begin{corollary}\label{Co:f bound}
For every $t>0$ and every function $v$ satisfying \eqref{E: Pois gen func cond},
\begin{align}
\nonumber&\exp\left \{-t\int_{\mathbb{R}^{2}}\bigl(1-v(x)\bigr)\Lambda(dx)\right \}\leq f_{\infty,\Lambda}(t,v) \\\noask
&\leq 1-\int_{\mathbb{R}^{2}}\dfrac{1-\exp\left \{-t\int_{\mathbb{R}^{2}}\Bigl(1-v(y)\bigl(1-h(x,y)\bigr)\Bigr)\Lambda(dy)\right \}}{\int_{\mathbb{R}^{2}}\Bigl(1-v(y)\bigl(1-h(x,y)\bigr)\Bigr)\Lambda(dy)}\bigl(1-v(x)\bigr)\Lambda(dx).
\end{align}
\end{corollary} 
\emph{Proof.}
The first inequality comes from the fact that $\mathcal{M}_{\infty}(\Phi_{t}) \subseteq \Phi_{t}$ a.s. and from Theorem \ref{T:PoissLaplace}. For the second inequality, we first use Theorem \ref{T:CSMA eq} to get
\begin{align*}
f_{\infty,\Lambda}(t,v)=1-\int_{0}^{t}\int_{\mathbb{R}^{2}}f_{\infty,\Lambda}\bigl(\tau,H(v,x)\bigr)\bigl(1-v(x)\bigr)\Lambda(dx)d\tau.
\end{align*}
 By the first inequality,
 \begin{align*}
 f_{\infty,\Lambda}\bigl(\tau,H(v,x)\bigr) \geq \exp\left \{-\tau\int_{\mathbb{R}^{2}}\Bigl(1-v(y)\bigl(1-h(x,y)\bigr)\Bigr)\Lambda(dy)\right \}.
 \end{align*}
Hence,
\begin{align*}
&f_{\infty,\Lambda}(t,v)\\\noask
&\leq 1-\int_{0}^{t}\int_{\mathbb{R}^{2}}\exp\left \{-\tau\int_{\mathbb{R}^{2}}\Bigl(1-v(y)\bigl(1-h(x,y)\bigr)\Bigr)\Lambda(dy)\right \}\bigl(1-v(x)\bigr)\Lambda(dx)d\tau.
\end{align*}
We then conclude by using Fubini's theorem. \eop

\subsection{$k$-Mat\'ern Models}\label{SS: k Matern}
Unlike the construction of the $\infty$-Mat\'ern model, which requires
reference to itself, the construction of the $k$-Mat\'ern model requires
reference to the $k-1$-Mat\'ern model. Thus, in order to study the p.g.fls
of the $k$-Mat\'ern model, we need to consider the joint distribution of 
all the $i$-Mat\'ern model for $i=0,\ldots,k$.
This can be done by considering the joint p.g.fls of these PPs:
\begin{align*}
f_{k,\Lambda}(t,v_{0},\ldots,v_{k})&:=G_{T_{t}(\mathcal{M}_{0}(\Phi)),\ldots,T_{t}(\mathcal{M}_{k}(\Phi))}(v_{0},\ldots,v_{k})\\\noask
&:=\textbf{E}\left [\prod_{i=0}^{k}\prod_{x \in T_{t}(\mathcal{M}_{i}(\Phi))}v_{i}(x)\right ].
\end{align*}
As the PPs $\mathcal{M}_{0}(\Phi)),\ldots,T_{t}(\mathcal{M}_{k}(\Phi))$ are not mutually disjoint, points in $T_{t}(\Phi)$ may appear more than once in the above product.Hence, computing this expression can be quite complicated. On the other hand, recall the following relation between $\mathcal{M}_{i}(\Phi)$,
\begin{itemize}
\item $\mathcal{M}_{0}(\Phi)\supseteq \mathcal{M}_{2}(\Phi)\cdots $;
\item $\mathcal{M}_{1}(\Phi)\supseteq \mathcal{M}_{3}(\Phi)\cdots $;
\item $\bigcup_{i=0}^{\infty}\mathcal{M}_{2i+1}(\Phi)\subseteq \bigcap_{i=0}^{\infty}\mathcal{M}_{2i}(\Phi)$.
\end{itemize}
We let
 \begin{align*}
\mathcal{Q}_{2i}(\Phi) &:=\mathcal{M}_{2i-2}(\Phi)\setminus \mathcal{M}_{2i}(\Phi);\\
  \mathcal{Q}_{2i-1}(\Phi) &:=\mathcal{M}_{2i-1}(\Phi)\setminus \mathcal{M}_{2i-3}(\Phi);\\ \mathcal{R}_{2i}(\Phi)&:=\mathcal{M}_{2i}(\Phi)\setminus \mathcal{M}_{2i-1}(\Phi);\\ 
\mathcal{R}_{2i-1}(\Phi)&:=\mathcal{M}_{2i-2}(\Phi)\setminus \mathcal{M}_{2i-1}(\Phi),
\end{align*}
 for $i=1,2,\ldots$ with the convention that $\mathcal{M}_{-1}(\Phi)=\emptyset$. We consider for each $k$ the functional
\begin{align}
\nonumber g_{k,\Lambda}(t,v_{1},\ldots,v_{k+1})&:=G_{T_{t}(\mathcal{Q}_{1}(\Phi)),\ldots,T_{t}(\mathcal{Q}_{k}(\Phi)),T_{t}(\mathcal{R}_{k}(\Phi))}\\\noask
&:=\mathbb{E}\left [\prod_{i=1}^{k}\prod_{x \in T_{t}(\mathcal{Q}_{i}(\Phi))}v_{i}(x)\prod_{x \in T_{t}(\mathcal{R}_{k}(\Phi))}v_{k}(x)\right ].
\end{align}
The next proposition shows that $f_{k,\Lambda}$ and $g_{k,\Lambda}$ are equivalent. Hence, we can study $g_{k,\Lambda}$ instead of $f_{k,\Lambda}$. As $\mathcal{Q}_{1}(\Phi),\ldots,\mathcal{Q}_{k}(\Phi),\mathcal{R}_{k}(\Phi)$ form a partition of $\Phi$, this makes $g_{k,\Lambda}$ easier to work with.
\begin{proposition}\label{P:f-g relation}
The functionals $f_{k,\Lambda}$ depend on the functionals $g_{k,\Lambda}$ in the sense that if we can compute $g_{k,\Lambda}(t,u_{1},\ldots,u_{k+1})$ for every $t \geq 0$ and for every $k+1$-tuple of functions $(u_{1},\ldots,u_{k+1})$, then we can compute $f_{k,\Lambda}(t,v_{0},\ldots,v_{k})$ for every $t \geq 0$ and for every $k+1$-tuple of functions $(v_{0},\ldots,v_{k})$.
\end{proposition}
\emph{Proof.}
Note that, by definition, $\mathcal{Q}_{i}(\Phi)$ must be either included in $\mathcal{M}_{j}(\Phi)$ or disjoint with $\mathcal{M}_{j}(\Phi)$ for every $i,j=0,1,\ldots$ Similarly, the same thing holds for $\mathcal{R}_{i}(\Phi)$ and $\mathcal{M}_{j}(\Phi)$. For each $i$, let $I_{i}:=\{j \st \mathcal{Q}_{i}(\Phi)\subseteq \mathcal{M}_{j}(\Phi) \}$ and $J_{k}:=\{j \st \mathcal{R}_{k}(\Phi)\subseteq \mathcal{M}_{j}(\Phi) \}$. In particular,
\begin{align*}
I_{2i-1}&=\{2i-1, 2i+1, \ldots \} \cup \{0, 2, 4,\ldots\}\\
I_{2i}&=\{0,2, \ldots 2i-2 \}\\
J_{2i-1} &=\{0,2, \ldots 2i-2 \}\\
J_{2i} &=\{0,2, \ldots 2i \}.
\end{align*}
 Then, fix a $k$, for every $i$, 
\begin{align*}
\mathcal{M}_{i}(\Phi)\left\{\begin{array}{ll}
                      =\displaystyle\left (\bigcup_{j \st i \in I_{j}}\mathcal{Q}_{j}(\Phi)\right )\cup\mathcal{R}_{k}(\Phi)  &\textnormal{ if } i \in J_{k} ,\\\noask
		      =\displaystyle\bigcup_{j \st i \in I_{j}}\mathcal{Q}_{j}(\Phi) &\textnormal{ otherwise}\,.
                     \end{array}\right.
\end{align*}
Hence,
\begin{align}
\nonumber &f_{k,\Lambda}(t,v_{0},\ldots,v_{k})\\\noask
\nonumber &=\textbf{E}\left [\prod_{i=0}^{k}\left (\prod_{x \in T_{t}(\mathcal{M}_{i}(\Phi))}v_{i}(x)\right )\right ]\\\noask
\nonumber &=\textbf{E}\left [\prod_{i=0}^{k}\left (\prod_{j \st i  \in I_{j}}\left (\prod_{x \in T_{t}(\mathcal{Q}_{j}(\Phi))}v_{i}(x)\right )\prod_{x \in  T_{t}(\mathcal{R}_{k}(\Phi))}\bigl(v_{i}(x)\bigr)^{\textbf{1}_{i \in J_{k}}}\right )\right  ]\\\noask
\nonumber &=\textbf{E}\left [\prod_{i=0}^{k}\left (\prod_{j =1}^{k}\left (\prod_{x \in T_{t}(\mathcal{Q}_{j}(\Phi))}\bigl(v_{i}(x)\bigr)^{\textbf{1}_{i \in I_{j}}}\right )\prod_{x \in T_{t}(\mathcal{R}_{k}(\Phi))}\bigl(v_{i}(x)\bigr)^{\textbf{1}_{i \in J_{k}}}\right )\right  ]\\\noask
&=\textbf{E}\left [ \prod_{j =1}^{k}\left (\prod_{x \in T_{t}(\mathcal{Q}_{j}(\Phi))} v_{i}(x)^{\sum_{i=0}^{k}\textbf{1}_{i \in I_{j}}}\right )\prod_{x \in T_{t}(\mathcal{R}_{k}(\Phi))}v_{i}(x)^{\sum_{i=0}^{k}\textbf{1}_{i \in J_{k}}}\right  ].
\end{align}
Put 
\begin{align*}
 u_{j}(x):=\prod_{i=0}^{k}v_{i}(x)^{\textbf{1}_{i \in I_{j}}},
\end{align*}
 for $j=1,\ldots, k$ and put 
 \begin{align*}
 u_{k+1}(x):=\prod_{i=0}^{k}v_{i}(x)^{\textbf{1}_{i \in J_{k}}}.
\end{align*}
The above equalities can be rewritten as
\begin{align*}
\nonumber &f_{k,\Lambda}(t,v_{0},\ldots,v_{k})=g_{k,\Lambda}(t,u_{1},\ldots,u_{k+1}).
\end{align*}
\eop

We now give a system of differential equations that governs the evolution of $g_{k,\Lambda}$ in $t$.
\begin{theorem}\label{T:k Matern eq}
For any locally finite measure $\Lambda$, the functional $g_{k,\Lambda}$ is continuous in $t$ and satisfies the following system of equations,
\begin{eqnarray}
\nonumber g_{k,\Lambda}(0,\textbf{v})&=&1;\\\noask
 \frac{dg_{k,\Lambda}(t,\textbf{v})}{dt}&=& -\sum_{i=1}^{k+1}\int_{\mathbb{R}^{2}}g_{k,\Lambda}\bigl(t,\textbf{H}_{i,k}(\textbf{v},x)\bigr)\bigl(w_{i}(x)-v_{i}(x)\bigr)\Lambda(dx),\label{E:k CSMA eq}
 \end{eqnarray}
 with $\textbf{v}=\{v_{1},\ldots,v_{k+1}\}$ is a $k+1$-tuple of functions from $\mathbb{R}^{2}$ to $[0,1]$ satisfying \eqref{E: gen def cond}, $(w_{1},\ldots,w_{k+1})$ is defined as follows
 \begin{align*}
 w_{i}&=v_{i+2} \mbox{ for $i$ odd smaller than $k+1$};\\
 w_{i}&=v_{i-2} \mbox{ for $i$ even smaller than $k+1$};\\
 w_{k+1}&=v_{k}
 \end{align*}
 when $k$ is even and
 \begin{align*}
 w_{i}&=v_{i+2} \mbox{ for $i$ odd smaller than $k$};\\
 w_{i}&=v_{i-2} \mbox{ for $i$ even smaller than $k+1$};\\
 w_{k}&=v_{k+1};\\
 w_{k+1}&=v_{k-1},
 \end{align*}
 when $k$ is odd where $v_{0}=\textbf{1}$ by convention and
 \begin{align}
 \textbf{H}_{i,k}(\textbf{v},x)(y)=(\mathfrak{u}_{0,i,k}(y,x),\cdots,\mathfrak{u}_{k+1,i,k}(y,x))
\end{align} 
where
\begin{align}
\nonumber \mathfrak{u}_{j,i,k}(y,x)&=v_{j}(y)\bigl(1-h(x,y)\bigr) \mbox{ if } i-1 \in I_{j} \mbox{ and } j \leq k;\\\noask
\nonumber \mathfrak{u}_{k+1,i,k}(y,x)&=v_{k+1}(y)\bigl(1-h(x,y)\bigr) \mbox{ if }  i-1 \in J_{k} ; \\\noask
\mathfrak{u}_{j,i,k}(y,x)&=v_{j}(y) \mbox{ otherwise,}
\end{align}
for every $y \in \mathbb{R}^{2}$.
\end{theorem}
\emph{Proof.} See Appendix \ref{A: Proof Theorem k Matern}.
\eop

\section{Reduced Palm Probability Generating Functional}\label{S: Palm int}
The aim if this section
is derive a similar system of differential
equations for the Palm versions of these
p.g.fls, i.e. the expectation of the product of
a function $v$ taken over all points of 
each of these point processes under its Palm distributions. 
We limit ourselves to the Mat\'ern models introduced in Section \ref{S:Model}.
We first derive the first moment measures of the Mat\'ern PPs from
their p.g.fls following a standard argument in the theory of PPs.
Then, by extending this method, the systems of integral
equations governing the evolution of the reduced Palm p.g.fls
of these PPs are also derived.

\subsection{Definitions}
 Let $\Phi$ be a PRRCR of intensity measure $\Lambda$ and expected conflict function $h$. We define
 \begin{align*}
& m_{k,\Lambda}(t,B)=\mathbb{E}\left [T_{t}\bigl(\mathcal{M}_{k}(\Phi)\bigr)(B)\right  ] \mbox{ for } k=\overline{0,\infty},
 \end{align*}
 for every bounded measurable set $B$. We can easily prove that for each $t>0$, $m_{k,\Lambda}(t,.),\, k=1,2,\ldots,\infty$ are indeed measures on $\mathbb{R}^{2}$ (see \cite{DaVer08}). They are called the first moment measures of the Mat\'ern models. As $T_{t}\bigl(\mathcal{M}_{k}(\Phi)\bigr)$ are thinning of $T_{t}(\Phi)$ for $k=\overline{0,\infty}$, their first moment measures must be absolutely continuous w.r.t. $\Lambda$. By abuse of notation , we denote by $m_{k,\Lambda}(t,x)$ the Radon-Nikodym derivative of $m_{k,\Lambda}$ w.r.t. $\Lambda$.
 
 Denote by $\textbf{P}_{x,\Phi}$ the reduced Palm distribution given a point at $x$ of a PP $\Phi$, respectively and denote by $\textbf{E}_{x,\Phi}$
the corresponding expectation. Note that we consider the reduced Palm distribution here rather than the Palm distribution (see Appendix \ref{App:AppendixA} for definitions).  Given the reduced Palm distribution, one can easily compute the Palm distribution by using Proposition \ref{P: Palm- reduced Palm}.

 We define
 \begin{align}
 f_{y,k,\Lambda}(t,v)=\textbf{E}_{y,T_{t}(\mathcal{M}_{k}(\Phi))}\left [\prod_{x \in T_{t}(\mathcal{M}_{k}(\Phi))} v(x)\right ]\mbox{ for } k=\overline{0,\infty}. 
 \end{align}
 Our objectives are to compute the measures $m_{k,\Lambda}$ and to derive systems of integral equations characterizing the functionals $f_{y,k,\Lambda}$.
 \subsection{First Moment Measure}
  These moment measures are computed in terms of the functionals $f_{k,\Lambda},\, k=\overline{0,\infty}$ as follows. 
 \begin{proposition}\label{P: First moment measure}
 For any locally finite measure $\Lambda$ and any $k$ from $1$ to $\infty$, we have,
 \begin{align}
 m_{k,\Lambda}(t,dx)=\int_{0}^{t}f_{k-1,\Lambda}\bigl(\tau,1-h(\cdot,x)\bigr).
 \end{align}  
 \end{proposition}
 \emph{Proof.}
 Since there is no ambiguity, we drop the subscript $\Lambda$ for notational convenience.
 
  We start with the case $k=\infty$. Applying Proposition \ref{P:Palm non Palm} with $v=e^{-s\textbf{1}_{B}}$, we get
 \begin{align*}
 \left .\frac{d}{ds}f_{\infty}(t,e^{-s\textbf{1}_{B}})\right |_{s=0}=-m_{\infty,\Lambda}(t,B).
 \end{align*}
  Moreover, by Theorem \ref{T:CSMA eq},
  \begin{align*}
  f_{\infty}(t,e^{-s\textbf{1}_{B}})=1-\int_{0}^{t}\int_{\mathbb{R}^{2}}f_{\infty,\Lambda}\bigl(\tau,H(e^{-s\textbf{1}_{B}},x)\bigr)(1-e^{-s\textbf{1}_{B}(x)})\Lambda(dx).
  \end{align*}
  Hence,
  \begin{align*}
m_{\infty}(t,B)=\left .\frac{d}{ds}\int_{0}^{t}\int_{\mathbb{R}^{2}}f_{\infty,\Lambda}\bigl(\tau,H(e^{-s\textbf{1}_{B}},x)\bigr)(1-e^{-s\textbf{1}_{B}(x)})\Lambda(dx)\right |_{s=0}.
  \end{align*}
Now we want to change the order of the derivative and the integration in the right hand side.
First, the conditions for this must be verified.
 
Using Proposition  \ref{P:Palm non Palm}, we have for all $s$,
\begin{align*}
 \left|\frac{d}{ds}f_{\infty}\bigl(\tau,H(e^{-s\textbf{1}_{B}},x)\bigr)\right|&=\left|\int_{B}f_{y,\infty}\bigl(\tau,H(e^{-s\textbf{1}_{B}},x)\bigr)m_{\infty}(\tau,y)\Lambda(dy)\right|\leq \tau\Lambda(B),
\end{align*}
since $m_{\infty,\Lambda}(\tau,x)<\tau$ for all $x$ (as the $\infty$-Mat\'ern model is a thinning of $\Phi$) and $0\leq H(e^{-s\textbf{1}_{B}},x)(y)\leq 1$ for every $x,y$ in $\mathbb{R}^{2}$. Hence,
\begin{align*}
&\left| \frac{d}{ds}\left(f_{\infty}\bigl(\tau,H(e^{-s\textbf{1}_{B}},x)\bigr)\left(1-e^{-s\textbf{1}_{B}(x)}\right)\right)\right|\\\noask
&\leq\textbf{1}_{B}(x)e^{-s\textbf{1}_{B}(x)}\left|f_{\infty}\bigl(\tau,H(e^{-s\textbf{1}_{B}},x)\bigr)\right|+\left(1-e^{-s\textbf{1}_{B}(x)}\right)\left|\frac{d}{ds}f_{\infty}\bigl(\tau,H(e^{-s\textbf{1}_{B}},x)\bigr)\right|\\\noask
&\leq \textbf{1}_{B}+s\textbf{1}_{B}\tau\Lambda(B) \leq \textbf{1}_{B}(x)\bigl(1+s\tau\Lambda(B)\bigr). 
\end{align*}
In the third line, we use the inequality $1-e^{-s\textbf{1}_{B}(x)} \leq s\textbf{1}_{B}(x)$. 
Since
\begin{align*}
&\int_{0}^{t}\int_{B}\bigl(1+s\tau\Lambda(B)\bigr)\,\Lambda(dx)\,d\tau= t\Lambda(B)+s\frac{t^{2}}{2}\Lambda(B)^{2} <\infty,
\end{align*}
we then have
\begin{align*}
&\int_{0}^{t}\int_{\mathbb{R}^{2}}\left| \frac{d}{ds}\left(f_{\infty}\bigl(\tau,H(e^{-s\textbf{1}_{B}},x)\bigr)\left(1-e^{-s\textbf{1}_{B}(x)}\right)\right)\right|\,\Lambda(dx)\,d\tau
< \infty.
\end{align*}
 So it is legitimate to change the order of
the differentiation and the integration,
\begin{align*}
\int_{B}m_{\infty}(t,x)\Lambda(dx)&=\int_{0}^{t}\int_{\mathbb{R}^{2}}\frac{d}{ds}f_{\infty}\bigl(\tau,H(e^{-s\textbf{1}_{B}},x)\bigr)\left(1-e^{-s\textbf{1}_{B}(x)}\right)\,\Lambda(dx)\,d\tau\Bigl|_{s=0}\\\noask
 &=\int_{0}^{t}\int_{B}f_{\infty}\bigl(\tau,H(e^{-s\textbf{1}_{B}},x)\bigr)e^{-s\textbf{1}_{B}(x)}\,\Lambda(dx)\,d\tau\Bigl|_{s=0}+\\\noask
&\hspace{.5cm}\int_{0}^{t}\int_{B}\frac{d}{ds}f_{\infty}\bigl(\tau,H(e^{-s\textbf{1}_{B}},x)\bigr)\left(1-e^{-s\textbf{1}_{B}(x)}\right)\,\Lambda(dx)\,d\tau\Bigl|_{s=0}\\\noask
 &=\int_{0}^{t}\int_{B}f_{\infty}\bigl(\tau,H(1,x)\bigr)\,\Lambda(dx)\,d\tau\\\noask
 &=\int_{B}\int_{0}^{t}f_{\infty}\bigl(\tau,H(1,x)\bigr)\,d\tau \,\Lambda(dx). \mbox{ (Fubini's theorem) }
\end{align*}
As this equality holds for any bounded measurable sets $B$, we must have
\[m_{\infty}(t,x)=\int_{0}^{t}f_{\infty}\bigl(\tau,H(1,x)\bigr)\,d\tau,\]
$\Lambda$-almost everywhere.

For $k<\infty$, by Proposition \ref{P:Palm non Palm},
\begin{align*}
\left .\frac{d}{ds}f_{k,}(t,e^{-s\textbf{1}_{B}})\right|_{s=0}=-m_{k}(t,B).
\end{align*} 
Moreover, by the same argument as in Proposition \ref{P:f-g relation},
 \[f_{k}(t,e^{-s\textbf{1}_{B}})= g_{k,}(t,\textbf{v}),\] where $\textbf{v}=(v_{1},\ldots,v_{k+1})$ with
 \begin{align*}
 v_{i}&=e^{-s\textbf{1}_{B}\textbf{1}_{k\in I(i)}} \mbox{ for } i=\overline{1,k};\\
 v_{k+1}&=e^{-s\textbf{1}_{B}\textbf{1}_{k\in J(i)}}.
\end{align*}

In particular, if $k$ is even, using the notation of Theorem \ref{T:k Matern eq},
\begin{align*}
&w_{i}=v_{i+2}=v_{i}=e^{-s\textbf{1}_{B}} \mbox{ for $i$ odd smaller than $k+1$};\\
&w_{i}=v_{i-2}=v_{i}=1 \mbox{ for $i$ even smaller than $k+1$};\\
&w_{k+1}=v_{k}=1;\\
&v_{k+1}=e^{-s\textbf{1}_{B}}.
\end{align*}
Hence,
\begin{align*}
\frac{d}{dt}f_{k}(t,v)=-\int_{\mathbb{R}^{2}}g_{k}\bigl(t,\textbf{H}_{k+1,k}(\textbf{v},x)\bigr)\bigl(1-e^{-s\textbf{1}_{B}}(x)\bigr)\Lambda(dx),
\end{align*}
which is equivalent to,
\begin{align*}
f_{k}(t,v)=1-\int_{0}^{t}\int_{\mathbb{R}^{2}}g_{k}\bigl(\tau,\textbf{H}_{k+1,k}(\textbf{v},x)\bigr)\bigl(1-e^{-s\textbf{1}_{B}}(x)\bigr)\Lambda(dx)d\tau.
\end{align*}
So,
\begin{align*}
m_{k}(t,B)=\left .\frac{d}{ds}\int_{0}^{t}\int_{\mathbb{R}^{2}}g_{k}\bigl(\tau,\textbf{H}_{k+1,k}(\textbf{v},x)\bigr)\bigl(1-e^{-s\textbf{1}_{B}}(x)\bigr)\Lambda(dx)d\tau\right|_{s=0}.
\end{align*}
By an argument similar to that of the case $k=\infty$, we can change the order of the derivative and the integrations  in the right hand side. As
\begin{align*}
&\left .\frac{d}{ds}g_{k}\bigl(\tau,\textbf{H}_{k+1,k}(\textbf{v},x)\bigr)\bigl(1-e^{-s\textbf{1}_{B}}(x)\bigr)\right|_{s=0}\\\noask
&=\left .\frac{d}{ds}g_{k}\bigl(\tau,\textbf{H}_{k+1,k}(\textbf{v},x)\bigr)\right|_{s=0}\bigl(1-e^{-s\textbf{1}_{B}}(0)\bigr)-\\
&\hspace{.5 cm}\left .g_{k}\bigl(\tau,\textbf{H}_{k+1,k}(\textbf{1},x)\bigr)\,\frac{d}{ds}\bigl(1-e^{-s\textbf{1}_{B}(x)}\bigr)\right|_{s=0}\\\noask
&=g_{k}\bigl(\tau,\textbf{H}_{k+1,k}(\textbf{1},x)\bigr)\textbf{1}_{B}(x),
\end{align*} 
we get
\begin{align*}
\int_{B}m_{k}(t,x)\Lambda(dx)=\int_{0}^{t}\int_{B}g_{k}\bigl(\tau,\textbf{H}_{k+1,k}(\textbf{1},x)\bigr)\Lambda(dx)d\tau.
\end{align*}
Note that $g_{k}\bigl(\tau,\textbf{H}_{k+1,k}(\textbf{1},x)\bigr)=f_{k-1}\Bigl(\tau,\bigl(1-h(.,x)\bigr)\Bigr)$. We get
\begin{align*}
\int_{B}m_{k}(t,x)\Lambda(dx)=\int_{B}\int_{0}^{t}f_{k-1}(\tau,(1-h(.,x)))d\tau\Lambda(dx),
\end{align*}
by Fubini's theorem. This proves the proposition for $k$ even. 
If $k$ is odd,
\begin{align*}
&w_{i}=v_{i+2}=v_{i}=v \mbox{ for $i$ odd smaller than $k$};\\
&w_{i}=v_{i-2}=v_{i}=1 \mbox{ for $i$ even smaller than $k+1$};\\
&w_{k+1}=v_{k-1}=v_{k-1}=1;\\
&w_{k}=v_{k+1}=1;\\
&v_{k}=v.
\end{align*}
By Theorem \ref{T:k Matern eq}, 
\begin{align*}
1-f_{k}(t,v)=-\int_{\mathbb{R}^{2}}g_{k}\bigl(t,\textbf{H}_{k,k}(\textbf{v},x)\bigr)\bigl(1-e^{-s\textbf{1}_{B}}\bigr)\Lambda(dx).
\end{align*}
The rest of the proof in this case is similar to the case $k$ is even, noting that $g_{k}(\tau,\textbf{H}_{k,k}\bigl(\textbf{1},x)\bigr)=f_{k-1}\Bigl(\tau,\bigl(1-h(.,x)\bigr)\Bigr)$ in this case.
\eop
\subsection{Reduced Palm Probability Generating Functionals}

We first derive a system of integral equations which have the Palm p.g.fls of the $\infty$-Mat\'ern model as a solution.
  \begin{proposition}\label{P:Palm gen infty Matern}
 For $\Lambda$-almost every $y$, for every positive functions $v$ satisfying (\ref{E: Pois gen func cond}) and every $t>0$, the functional $f_{y,\Lambda}(t,v)$ satisfies the system of integral equations
\begin{align*}
& f_{y,\infty,\Lambda}(0,v)=1;\\\noask
& f_{y,\infty,\Lambda}(t,v)=\frac{\int_{0}^{t}f_{\Lambda}\bigl(\tau,H(v,y)\bigr)d\tau}{m_{\infty,\Lambda}(t,y)}\\
&-\int_{0}^{t}\int_{\mathbb{R}^{2}}f_{y,\infty,\Lambda}\bigl(t,H(v,x)\bigr)\bigl(1-v(x)\bigr)\bigl(1-h(x,y)\bigr)\frac{m_{\infty,\Lambda}(\tau,y)}{m_{\infty,\Lambda}(t,y)}\Lambda(dx)d\tau,
\end{align*}
where $H$ is defined in (\ref{E:h transform}) and $m_{\infty,\Lambda}(t,y)$ is the Radon-Nikodym derivative w.r.t. $\Lambda$ of the intensity measure of 
$T_{t}\bigl(\mathcal{M}_{\infty}(\Phi)\bigr)$ computed in Proposition \ref{P: First moment measure}.
\end{proposition}
 \emph{Proof.} See Appendix \ref{A: Proof Proposition Palm}.
 
As the study of the $k$-Mat\'ern models p.g.fls requires the joint p.g.fls of the PPs induced by the $0$-Mat\'ern, $1$-Mat\'ern, $\ldots$, $k$-Mat\'ern models, it is natural that the study of its Palm version requires the \emph{Palm version} of theses joint p.g.fls. Defining such Palm joint p.g.fls is quite complicated in general. Nevertheless, to exemplify the idea, we present here the result for $k=\infty$ and $k=1$.

First of all, we can assume w.l.o.g. that $u(x)>0$ for any $x$ in $\mathbb{R}^{2}$ (as $G_{T_{t}(\mathcal{M}_{1}(\Phi))}(v)$ is obtained by setting $u=\textbf{1}$ in $g_{1,\Lambda}(t,v,u)$ and such $u$ satisfies this condition). By the same arguments as in Proposition \ref{P:Palm non Palm},
\begin{align*}
&\left .-\frac{d}{ds}g_{1,\Lambda}(t,ve^{-s\textbf{1}_{B}},u)\right |_{s=0}\\\noask
&=\mathbb{E}\left [ T_{t}\bigl(\mathcal{Q}_{1}(\Phi)\bigr)(B)\prod_{y \in  T_{t}(\mathcal{Q}_{1}(\Phi))}v(y)\prod_{y\in T_{t}(\mathcal{R}_{1}(\Phi))}u(y)\right ]\\\noask
&=\mathbb{E}\left [ T_{t}\bigl(\mathcal{M}_{1}(\Phi)\bigr)(B)\prod_{y \in  T_{t}(\mathcal{M}_{1}(\Phi))}v(y)\prod_{y\in T_{t}(\Phi)\setminus T_{t}(\mathcal{M}_{1}(\Phi))}u(y)\right ]\\\noask
&=\mathbb{E}\left [\sum_{x \in \mathcal{M}_{1}(\Phi)}\textbf{1}_{x\in B}\prod_{y \in  \mathcal{M}_{1}(\Phi_{t})}v(y)\prod_{y\in \Phi_{t}}u(y)\right ]\\\noask
&=\mathbb{E}\left [\sum_{x \in T_{t}(\mathcal{M}_{1}(\Phi))\cap B}v(x)\prod_{y \in  T_{t}(\mathcal{M}_{1}(\Phi))\setminus \{x\}}v(y)\prod_{y\in \Phi_{t}}u(y)\right ]\\\noask
&\leq \mathbb{E}\left [\sum_{x \in \mathcal{M}_{1}(\Phi_{t})\cap B}v(x)\right ]= \int_{B}v(x)m_{1,\Lambda}(t,x)\Lambda(dx).
\end{align*}
 The second equality comes from the fact that $\mathcal{Q}_{1}(\Phi)=\mathcal{M}_{1}(\Phi)$ and $\mathcal{R}_{1}(\Phi)=\Phi \setminus \mathcal{M}_{1}(\Phi)$. 
So, for each $u$, $v$, $\left .-\frac{d}{ds}g_{\Lambda}(t,ve^{-s\textbf{1}_{.}},u)\right |_{s=0}$, considered as a measure in $\mathbb{R}^{2}$, is absolutely continuous w.r.t. $v(x)m_{1,\Lambda}(t,x)\Lambda(dx)$. Hence, it admits a Radon-Nikodym derivative w.r.t. the latter, i.e.
 \begin{align}
\left .-\frac{d}{ds}g_{\Lambda}(t,ve^{-s\textbf{1}_{B}},u)\right |_{s=0}=\int_{B}g_{x,1,\Lambda}(t,v,u)v(x)m_{1,\Lambda}(t,x)\Lambda(dx).\label{E:join Palm-non Palm prod}
\end{align} 
Now, by the same arguments as in Proposition \ref{P:reduced Palm non Palm} , we can also prove that for every measurable set $B$ not in the support of $v$, 
\begin{align}
\left .\frac{d}{ds}g_{\Lambda}(t,v+s\textbf{1}_{B},u)\right |_{s=0}=\int_{B}g_{x,1,\Lambda}(t,v,u)m_{1,\Lambda}(t,x)\Lambda(dx).\label{E:join Palm-non Palm sum}
\end{align}

In particular, by taking $u=\textbf{1}$, we get:
\begin{proposition}
For every $x$ every function $v$ satisfying \eqref{E: Pois gen func cond}, 
\begin{align}
G^{!}_{x,\mathcal{M}_{1}(\Phi_{t})}(v)=g_{x,1,\Lambda}(t,v,\textbf{1}).
\end{align}
\end{proposition}
\emph{Proof.}
By Proposition \ref{P:Palm non Palm} and Proposition \ref{P:reduced Palm non Palm},
\begin{align*}
-\left .\frac{d}{ds}g(t,ve^{-s\textbf{1}_{B}},\textbf{1})\right |_{s=0}&=-\left .\frac{d}{ds}G_{T_{t}(\mathcal{M}_{1}(\Phi))}(ve^{-s\textbf{1}_{B}})\right |_{s=0}\\\noask
&=\int_{B}v(x)G^{!}_{x,\mathcal{M}_{1}(\Phi_{t})}(v)m_{1,\Lambda}(t,x)\,\Lambda(dx)\\\noask
\left .\frac{d}{ds}g(t,v+s\textbf{1}_{B},\textbf{1})\right |_{s=0}&=\left .\frac{d}{ds}G_{T_{t}(\mathcal{M}_{1}(\Phi))}(v+s\textbf{1}_{B})\right |_{s=0}\\\noask
&=\int_{B}G^{!}_{x,\mathcal{M}_{1}(\Phi_{t})}(v)m_{1,\Lambda}(t,x)\,\Lambda(dx).
\end{align*}
By comparing these equalities with \eqref{E:join Palm-non Palm prod} and \eqref{E:join Palm-non Palm sum}, we can conclude the proof. \eop

Note that the equations \eqref{E:join Palm-non Palm prod} and \eqref{E:join Palm-non Palm sum} are very similar with the result of Proposition \ref{P:Palm non Palm}, hence we can call $g_{x,\Lambda}$ the \emph{Palm version} of the joint p.g.fl $g_{1,\Lambda}$ by abuse of notation. 

We next derive a system of equations that has $g_{x,1,\Lambda}(t,u,v)$ as a solution in the same spirit as in Proposition \ref{P:Palm gen infty Matern}.

\begin{proposition}\label{P: type II reduced Palm}
For $\Lambda$-almost every $y$, every functions $u>0$, $v$ satisfying \eqref{E: Pois gen func cond}, the functional $g_{y,\Lambda}$ satisfies the integral equation
\begin{align*}
&g_{y,\Lambda}(0,v,u)=1;\\\noask
&g_{y,\Lambda}(t,v,u)= \int_{0}^{t}\frac{g_{\Lambda}(\tau,H(v,y),H(u,y))}{m_{1,\Lambda}(t,y)}d\tau-\int_{0}^{t}\int_{\mathbb{R}^{2}}\Bigl (g_{y,\Lambda}(\tau,v,u)\bigl(1-u(x)\bigr)+\\\noask
&g_{y,\Lambda}(\tau,H(v,x),H(u,x))\bigl(u(x)-v(x)\bigr)\bigl(1-h(x,y)\bigr)\Bigr )\frac{m_{1,\Lambda}(\tau,y)}{m_{1,\Lambda}(t,y)}\Lambda(dx)d\tau,
\end{align*}
where $H$ is defined in (\ref{E:h transform}).
\end{proposition} 
\emph{Proof.} We present here only a sketch of this proof. To complete this sketch, it is sufficient to prove the condition for changing the order of derivatives and integrations in the same way as in the proof of Proposition \ref{P:Palm gen infty Matern}. As there is no ambiguity, we drop the $1$ and $\Lambda$ subscript.
 By definition, for any bounded measurable $B$ in the support of $v$,
\begin{align*}
\frac{d}{ds}g(t,ve^{-s\textbf{1}_{B}},u)\Biggr|_{s=0}=-\int_{B}g_{x}(t,v,u)v(x)m(t,x)\Lambda(dx).
\end{align*}
By Theorem \ref{T:k Matern eq},
\begin{align*}
g(t,ve^{-s\textbf{1}_{B}},u)=&1-\int_{0}^{t}\int_{\mathbb{R}^{2}}g\left (\tau,ve^{-s\textbf{1}_{B}},u\right )\bigl(1-u(x)\bigr)\Lambda(dx)d\tau-\\\noask
&\int_{0}^{t}\int_{\mathbb{R}^{2}}g\left (\tau,H(ve^{-s\textbf{1}_{B}},x),H(u,x)\right )\left (u(x)-v(x)e^{-s\textbf{1}_{B}(x)}\right )\Lambda(dx)d\tau.
\end{align*}
So,
\begin{align*}
&\int_{B}g_{x}(t,v,u)v(x)m(t,x)\Lambda(dx)\\\noask
&=\frac{d}{ds}\int_{0}^{t}\int_{\mathbb{R}^{2}}g_{\Lambda}\left (\tau,ve^{-s\textbf{1}_{B}},u\right )\bigl(1-u(x)\bigr)\Lambda(dx)d\tau\Biggr|_{s=0}\\\noask
&\;\;\;+\frac{d}{ds}\int_{0}^{t}\int_{\mathbb{R}^{2}}g\left (\tau,H(ve^{-s\textbf{1}_{B}},x),H(u,x)\right )\left (u(x)-v(x)e^{-s\textbf{1}_{B}(x)}\right )\Lambda(dx)d\tau\Biggr|_{s=0}.
\end{align*}
The first term is
\begin{align}
\nonumber &\frac{d}{ds}\int_{0}^{t}\int_{\mathbb{R}^{2}}g_{\Lambda}\left (\tau,ve^{-s\textbf{1}_{B}},u\right )\bigl(1-u(x)\bigr)\Lambda(dx)d\tau\Biggr|_{s=0}\\\noask
\nonumber&=\int_{0}^{t}\int_{\mathbb{R}^{2}}\frac{d}{ds}g\left (\tau,ve^{-s\textbf{1}_{B}},u\right )\Biggr|_{s=0}\bigl(1-u(x)\bigr)\Lambda(dx)d\tau\\\noask
\nonumber&=-\int_{0}^{t}\int_{\mathbb{R}^{2}}\left (\int_{B}g_{y}(\tau,v,u)v(y)m(\tau,y)\Lambda(dy)\right )\bigl(1-u(x)\bigr)\Lambda(dx)d\tau\\\noask
&=-\int_{B}\left (\int_{0}^{t}\int_{\mathbb{R}^{2}}g_{y}(\tau,v,u)\bigl(1-u(x)\bigr)\Lambda(dx)d\tau\right )v(y)m(\tau,y)\Lambda(dy).\label{E: reduced g e1}
\end{align}
For the second term,
\begin{align*}
&\frac{d}{ds}\int_{0}^{t}\int_{\mathbb{R}^{2}}g\left (\tau,H(ve^{-s\textbf{1}_{B}},x),H(u,x)\right )\left (u(x)-v(x)e^{-s\textbf{1}_{B}(x)}\right )\Lambda(dx)d\tau\Biggr|_{s=0}\\\noask
&+\int_{0}^{t}\int_{\mathbb{R}^{2}}\frac{d}{ds}g\left (\tau,H(ve^{-s\textbf{1}_{B}},x),H(u,x)\right )\left (u(x)-v(x)e^{-s\textbf{1}_{B}(x)}\right )\biggr|_{s=0}\Lambda(dx)d\tau\\\noask
&=\int_{0}^{t}\int_{\mathbb{R}^{2}}\frac{d}{ds}g\left (\tau,H(ve^{-s\textbf{1}_{B}},x),H(u,x)\right )\biggr|_{s=0}\left (u(x)-v(x)e^{-s\textbf{1}_{B}(x)}\right )\biggr|_{s=0}\Lambda(dx)d\tau\\\noask
&\;\;\;+\int_{0}^{t}\int_{\mathbb{R}^{2}}g\left (\tau,H(ve^{-s\textbf{1}_{B}},x),H(u,x)\right )\biggr|_{s=0}\frac{d}{ds}\left (u(x)-v(x)e^{-s\textbf{1}_{B}(x)}\right )\biggr|_{s=0}\Lambda(dx)d\tau.
\end{align*}
The second term in the above equality is
\begin{align}
\nonumber&\int_{0}^{t}\int_{\mathbb{R}^{2}}g\left (\tau,H(ve^{-s\textbf{1}_{B}},x),H(u,x)\right )\biggr|_{s=0}\frac{d}{ds}\left (u(x)-v(x)e^{-s\textbf{1}_{B}(x)}\right )\biggr|_{s=0}\Lambda(dx)d\tau\\\noask
\nonumber&=\int_{0}^{t}\int_{\mathbb{R}^{2}}g(\tau,v,H(u,x))v(x)\textbf{1}_{B}(x)\Lambda(dx)d\tau\\\noask
&=\int_{B}\int_{0}^{t}g(\tau,H(v,x),H(u,y))v(y)d\tau\Lambda(dy), \label{E: reduced g e2}
\end{align}
while the first term is
\begin{align}
\nonumber& \int_{0}^{t}\int_{\mathbb{R}^{2}}\frac{d}{ds}g\left (\tau,H(ve^{-s\textbf{1}_{B}},x),H(u,x)\right )\biggr|_{s=0}\left (u(x)-v(x)e^{-s\textbf{1}_{B}(x)}\right )\biggr|_{s=0}\Lambda(dx)d\tau\\\noask
\nonumber& =-\int_{0}^{t}\int_{\mathbb{R}^{2}}\int_{B}g_{y}\bigl(\tau,H(v,x),H(u,x)\bigr)H(v,x)(y)m(\tau,y)\Lambda(dy)\bigl(u(x)-v(x)\bigr)\\
 &\hfill \nonumber \Lambda(dx)d\tau\\\noask 
\nonumber& =-\int_{0}^{t}\int_{\mathbb{R}^{2}}\int_{B}g_{y}\bigl(\tau,H(v,x),H(u,x)\bigr)\bigl(1-h(x,y)\bigr)v(y)m(\tau,y)\Lambda(dy)\\
&\hfill \nonumber \bigl(u(x)-v(x)\bigr)\Lambda(dx)d\tau\\\noask
\nonumber& =-\int_{B}\left (\int_{0}^{t}\int_{\mathbb{R}^{2}}g_{y}\bigl(\tau,H(v,x),H(u,x)\bigr)\bigl(u(x)-v(x)\bigr)\bigl(1-h(x,y)\bigr)\Lambda(dx)d\tau\right )\\
&  v(y)m(\tau,y)\Lambda(dy). \label{E: reduced g e3}
\end{align}
Putting together \eqref{E: reduced g e1}, \eqref{E: reduced g e2} and \eqref{E: reduced g e3}, we get
\begin{align*}
&\int_{B}g_{y}(t,v,u)v(y)m(t,y)\Lambda(dy)=\int_{B}\int_{0}^{t}g(\tau,H(v,x),H(u,y))v(y)d\tau\Lambda(dy)\\\noask
&-\int_{B}\left (\int_{0}^{t}\int_{\mathbb{R}^{2}}g_{y}(\tau,v,u)\bigl(1-u(x)\bigr)\Lambda(dx)d\tau\right )v(y)m(\tau,y)\Lambda(dy)-\int_{B}\left (\int_{0}^{t}\right .\\\noask
&\int_{\mathbb{R}^{2}}g_{y}(\tau,H(v,x),H(u,x))\bigl(u(x)-v(x)\bigr)\bigl(1-h(x,y)\bigr)\Lambda(dx)d\tau\biggr )v(y)m(\tau,y)\Lambda(dy)
\end{align*}
for every bounded measurable $B$ in the support of $v$, which proves the proposition for $\Lambda$-almost every $x$ such that $v(x)>0$. For $x$ such that $v(x)=0$, proceed as above and use the equality
\begin{align*}
\left .\frac{d}{ds}g(t,v+s\textbf{1}_{B},u)\right |_{s=0}=\int_{B}g_{x}(t,v,u)m(t,x)\Lambda(dx).
\end{align*}
for every bounded measurable $B$ not in the support of $v$. \eop

\emph{Remarks:} Note that the method used in this section can be extended to compute the $n^{th}$-moment measures and the $n^{th}$-fold Palm distribution of the PPs induced by the Mat\'ern models. In particular, the $n^{th}$ moment measure is the $n^{th}$ derivative of the p.g.fls at $0$ for a suitably chosen function $v$. Thus, by using Theorems \ref{T:CSMA eq} and  \ref{T:k Matern eq}, we can compute these measures by an induction argument on $n$. Then, by using the same argument but with a general function $v$, we can derive a system of integral equations that has the $n^{th}$-fold Palm p.g.fls as solutions. These computations are, however, long and tedious and are hence omitted here.

\section*{Acknowledgment}
The authors would like to thank D. Daley for his
very valuable comments on this work.
The early stages of this work were initiated 
in France at Ecole Normale Sup\'erieure and INRIA,
The later stages were pursued in the USA
at Qualcomm for the first author and 
at the University of Texas at Austin for the second author,
where they were supported
by a grant of the Simons Foundation (\#197982 to UT Austin). 
 
 \bibliographystyle{acm}
\bibliography{Bib}  
\newpage
\appendix
\section{Palm version of Probability Generating Functionals}\label{A: Palm}

 We define the Palm version and the reduced Palm version of the p.g.fls of a PP $N$ as
 \begin{align}
 G_{x,N}(v)&:=\mathbb{E}_{x,N}\left [\prod_{y\in N}v(y)\right ]\label{E:reduced Palm pgfl};\\\noask
  G^{!}_{x,N}(v)&:=\mathbb{E}^{!}_{x,N}\left [\prod_{y\in N}v(y)\right ]\label{E:Palm pgfl},
 \end{align}
 where $\mathbb{E}_{x,N}$ and $\mathbb{E}^{!}_{x,N}$ are the Palm distribution and the reduced Palm distribution of $N$, respectively.
The relation between the reduced Palm version and the Palm version is characterized as follows.
\begin{proposition}\label{P: Palm- reduced Palm}
For any PP $N$ with intensity measure $m_{N}$,we have for $m_{N}$-almost every $x$,
\begin{align}
  G_{x,N}(v)=v(x)G^{!}_{x,N}(v).
\end{align}
\end{proposition}
\emph{Proof}
Let $f(x,N)= \prod_{y \in N}v(y)$ and $g(x,N)= v(x)\prod_{y \in N}v(y)$. We have for any $x \in N$,
\[f(x,N)=v(x)f(x,N\setminus x)=g(x,N\setminus x).\]
So,
\begin{align*}
\int_{\mathbb{R}^{2}}\mathbb{E}_{x}\left [f(x,N)\right ]m_{N}(dx)&= \mathbb{E}\left [\sum_{x\in N}f(x,N)\right ]\\\noask
    &= \mathbb{E}\left [\sum_{x\in N}g(x,N \setminus x)\right ]\\\noask
    &= \int_{\mathbb{R}^{d}}\mathbb{E}^{!}_{x}\left [g(x,N)\right ]m_{N}(dx)\\\noask
    &= \int_{\mathbb{R}^{d}}\mathbb{E}^{!}_{x}\left [v(x)f(x,N)\right ]m_{N}(dx)\\\noask
    &= \int_{\mathbb{R}^{d}}v(x)\mathbb{E}^{!}_{x}\left [f(x,N)\right ]m_{N}(dx).
\end{align*}
The conclusion then follows directly.
\eop 

Hence, we can write
\begin{align}
G^{!}_{x,N}(v)&:=\mathbb{E}^{!}_{x,N}\left [\prod_{y\in N}v(y)\right ]=\mathbb{E}^{!}_{x,N}\left [\prod_{y\in N\setminus\{x\}}v(y)\right ].
\end{align} 
And moreover, it is now sufficient to concentrate on the reduced Palm p.g.fls.
The next two results give us the relation between the Palm and non Palm versions of the p.g.fl of a PP.
\begin{proposition}\label{P:Palm non Palm}
Let $N$ be a PP with locally finite intensity measure. Then for any function $v$ and any bounded measurable set $B$,
\begin{align}
\nonumber \frac{d}{dt}G_{N}\left (ve^{-t\textbf{1}_{B}}\right )&=-\int_{B}G_{x,N}\left (ve^{-t\textbf{1}_{B}}\right )m_{N}(dx)\\\noask
&=-\int_{B}v(x)e^{-t\textbf{1}_{ B}(x)}G^{!}_{x,N}\left (ve^{-t\textbf{1}_{B}}\right )m_{N}(dx).
\end{align}
\end{proposition}
\emph{Proof.}
We have,
\begin{align*}
\frac{d}{dt}\left (\prod_{x \in N}v(x)e^{-t\textbf{1}_{B}(x)}\right )&= \frac{d}{dt}\left (\prod_{x \in N}v(x)\right )e^{-t N(B)}\\\noask
&=-\left (\prod_{x \in N}v(x)\right )e^{-t N(B)}N(B)\\\noask
&=-\left (\prod_{x \in N}v(x)e^{-t\textbf{1}_{B}(x)}\right )N(B)\\\noask
&= -\sum_{x \in B \cap N}\left(\prod_{y \in N}v(y)e^{-t\textbf{1}_{B}(y)}\right ).
\end{align*}
Since the absolute value of the last line is bounded above by $N(B)$ and $\mathbb{E}[N(B)]= m_{N}(B)<\infty$ by definition, we have by bounded convergence theorem,
\begin{align*}
\frac{d}{dt}G_{N}(ve^{-t\textbf{1}_{B}})&=\frac{d}{dt}\mathbb{E}\left [\left (\prod_{x \in N}v(x)e^{-t\textbf{1}_{B}(x)}\right )\right ]\\\noask
&=\mathbb{E}\left [\frac{d}{dt}\left (\prod_{x \in N}v(x)e^{-t\textbf{1}_{B}(x)}\right )\right]\\\noask
&=-\mathbb{E}\left [\sum_{x \in B \cap N}\left(\prod_{y \in N}v(y)e^{-t\textbf{1}_{B}(y)}\right )\right].
\end{align*}
As the last line is equal to $-\int_{B}G_{x,N}(ve^{-t\textbf{1}_{B}})m_{N}(dx)$ by the Campbell formula, the conclusion follows directly. 
\eop

However, the previous proposition does not give us the relation between the p.g.fl and its reduced Palm versions at $x$ such that $v(x)=0$. For such $x$, we need the following result.
\begin{proposition}\label{P:reduced Palm non Palm}
Let $N$ be a PP with locally finite intensity measure and $t$ be positive number smaller than $1$.  For any function $v$ and any measurable set $B$ not in the support of $v$ (i.e. $v(x)=0$ for every $x$ in $B$),
\begin{align}
\frac{d}{dt}G_{N}(v+t\textbf{1}_{B})=\int_{B}G^{!}_{x,N}(v+t\textbf{1}_{B})m_{N}(dx).
\end{align}
\end{proposition}  
\emph{Proof.}
Since $B$ is not in the support of $v$,
\begin{align*}
\frac{d}{dt}\left (\prod_{y\in N}\bigl(v(y)+\textbf{1}_{B}(y)\bigr)\right )&=\frac{d}{dt}\left (\prod_{y\in N \setminus B}v(y)\right )t^{N(B)}\\\noask
&=\left (\prod_{y\in N \setminus B}v(y)\right )t^{N(B)-1}N(B)\\\noask
&=\left (\prod_{y\in N \setminus B}v(y)\right )\left (\sum_{x\in N \cap B}\prod_{y\in (N \cap B) \setminus \{x\}}t\right )\\\noask
&=\sum_{x\in N \cap B}\left (\prod_{y\in N \setminus B}v(y)\prod_{y\in N \cap B \setminus \{x\}}t\right )\\\noask
&=\sum_{x\in N \cap B}\left (\prod_{y\in N \setminus \{x\}}\bigl(v(y)+\textbf{1}_{B}(y)\bigr)\right ).
\end{align*}
Again, by bounded convergence theorem,
\begin{align*}
\frac{d}{dt}G_{N}(v+t\textbf{1}_{B})&=\mathbb{E}\left [\frac{d}{dt}\left (\prod_{y\in N}\bigl(v(y)+\textbf{1}_{B}\bigr)\right )\right]\\\noask
&=\mathbb{E}\left [\sum_{x\in N \cap B}\left (\prod_{y\in N \setminus x}\bigl(v(y)+\textbf{1}_{B}\bigr)\right )\right].
\end{align*}
The conclusion follows directly from the fact that the last line equals \[\int_{B}G^{!}_{x,N}(v+t\textbf{1}_{B})m_{N}(dx)\] 
by reduced Campbell formula.
\eop

\section{Proof of Lemma \ref{L:limit 0}}\label{A: Proof Lemma Limit 0}
Every expectation in this proof should be understood as the conditional expectation given $\Phi_{t}$. 
By Lemma \ref{L: main eq},
\begin{align*}
 \prod_{x \in T_{t,t+\epsilon}(\mathcal{M}_{\infty}(\Phi))}v(x)=1+\sum_{i=1}^{\infty}(-1)^{i}\sum_{(x_{1},\ldots,x_{i}) \in \left(T_{t,t+\epsilon}(\mathcal{M}_{\infty}(\Phi))\right)^{(i!)}}\prod_{j=1}^{i}\bigl(1-v(x_{j})\bigr).
\end{align*}
The first step is to show that for $\epsilon$ small enough,
\begin{align*}
 \left|\mathbb{E}\left[\sum_{i=2}^{\infty}\sum_{(x_{1},\ldots,x_{i}) \in \left(T_{t,t+\epsilon}(\mathcal{M}_{\infty}(\Phi))\right)^{(i!)}}\prod_{j=1}^{i}\bigl(1-v(x_{j}) \right]\right| \leq 2\epsilon^{2}\left(\int_{\mathbb{R}^{2}}\bigl(1-v(x)\bigr)\,\Lambda(dx)\right)^{2}\\ \mbox{ a.s. }
\end{align*}
 Since $T_{t,t+\epsilon}\bigl(\mathcal{M}_{\infty}(\Phi)\bigr) \subseteq T_{t,t+\epsilon}(\Phi)$ a.s., 
\begin{align*}
& \left|\mathbb{E}\left[\sum_{i=2}^{\infty}\sum_{(x_{1},\ldots,x_{i}) \in (T_{t,t+\epsilon}(\mathcal{M}_{\infty}(\Phi)))^{(i!)}}\prod_{j=1}^{i}\bigl(1-v(x_{j})\bigr)\right]\right|\\\noask
&\leq\; \mathbb{E}\left[\sum_{i=2}^{\infty}\sum_{(x_{1},\ldots,x_{i}) \in (T_{t,t+\epsilon}(\Phi))^{(i!)}}\prod_{j=1}^{i}\bigl(1-v(x_{j})\bigr)\right]\\\noask
&=\; \sum_{i=2}^{\infty}\mathbb{E}\left[\sum_{(x_{1},\ldots,x_{i}) \in (T_{t,t+\epsilon}(\Phi))^{(i!)}}\prod_{j=1}^{i}\bigl(1-v(x_{j})\bigr)                                                                                   \right]\\\noask
&=\; \sum_{i=2}^{\infty}\epsilon^{i}\left(\int_{\mathbb{R}^{2}}\bigl(1-v(x)\bigr)\,\Lambda(dx)\right)^{i} \mbox{ a.s.},
\end{align*}
where in the last line we apply the multivariate Campbell formula to the PPP $T_{t,t+\epsilon}(\Phi)$ \cite[p.112]{SKM95}. Take now any $\epsilon < \frac{1}{2}\left(\int_{\mathbb{R}^{2}}\bigl(1-v(x)\bigr)\,\Lambda(dx)\right)^{-1}$,  
\begin{align*}
 \sum_{i=2}^{\infty}\epsilon^{i}\left(\int_{\mathbb{R}^{2}}\bigl(1-v(x)\bigr)\,\Lambda(dx)\right)^{i}&=\epsilon^{2}\frac{\left(\int_{\mathbb{R}^{2}}\bigl(1-v(x)\bigr)\,\Lambda(dx)\right)^{2}}{1-\epsilon\int_{\mathbb{R}^{2}}\bigl(1-v(x)\bigr)\,\Lambda(dx)}\\\noask
& \leq 2\epsilon^{2}\left(\int_{\mathbb{R}^{2}}\bigl(1-v(x)\bigr)\,\Lambda(dx)\right)^{2}.
\end{align*}

The next step is to bound $\mathbb{E}\left[\sum_{x \in T_{t,t+\epsilon}\bigl(\mathcal{M}_{\infty}(\Phi)\bigr)}\bigl(1-v(x)\bigr)\right]$.
 By Proposition \ref{P:inf thin layer},
\begin{align*}
\Delta_{\infty,d,t,t+\epsilon} \subseteq T_{t,t+\epsilon}\bigl(\mathcal{M}_{\infty}(\Phi)\bigr) \subseteq \Delta_{\infty,u,t,t+\epsilon} \mbox{ a.s.}
\end{align*}
Then,
\begin{align*}
& \sum_{x \in \Delta_{\infty,d,t,t+\epsilon}}\bigl(1-v(x)\bigr)\leq \sum_{x \in T_{t,t+\epsilon}\bigl(\mathcal{M}_{\infty}(\Phi)\bigr)}\bigl(1-v(x)\bigr)\leq \sum_{x \in \Delta_{\infty,u,t,t+\epsilon}}\bigl(1-v(x)\bigr)\mbox{ a.s.}
\end{align*}
Given $\Phi_{t}$, $\Delta_{\infty,u,t,t+\epsilon}$ is an independently thinning of $T_{t,t+\epsilon}(\Phi)$ with thinning probability $\prod_{y \in T_{t}(\mathcal{M}_{\infty}(\Phi))}\bigl(1-h(x,y)\bigr)$. Hence, it is a PPP of intensity $\epsilon\,\Lambda_{\infty,t}$. We have then,
\begin{align}
\mathbb{E}\left[\sum_{x \in \Delta_{\infty,u,t,t+\epsilon}}\bigl(1-v(x)\bigr)\right]=\epsilon\int_{\mathbb{R}^{2}}\bigl(1-v(x)\bigr)\,\Lambda_{\infty,t}(dx) \mbox{ a.s. } \label{E: thin upper bound}
\end{align}
 Moreover, we can compute the intensity of $ \Delta_{\infty,d,t,t+\epsilon}$ (conditioned on $\Phi_{t}$) as follow. Take any bounded measurable set $A$ in $\mathbb{R}^{2}$,
\begin{align*}
& \mathbb{E}\left[\left| \Delta_{\infty,d,t,t+\epsilon}\cap A\right|\right]\\\noask
&= \mathbb{E}\left[\sum_{x \in T_{t,t+\epsilon}(\Phi)\cap A}\left(\left (\prod_{y \in T_{t,t+\epsilon}(\Phi)}\textbf{1}_{C(x,y)=0}\right )\left (\prod_{y \in T_{t}(\mathcal{M}_{\infty}(\Phi))}\textbf{1}_{C(x,y)=0}\right )\right )\right]\\\noask
&= \mathbb{E}\left[\sum_{x \in T_{t,t+\epsilon}(\Phi)\cap A}\left(\left (\prod_{y \in T_{t,t+\epsilon}(\Phi)}\textbf{1}_{C(x,y)=0}\right )\left (\prod_{y \in T_{t}(\mathcal{M}_{\infty}(\Phi))}\bigl(1-h(x,y)\bigr)\right )\right )\right].
\end{align*}
Let $\mathbb{P}_{x}^{!}$ be the \emph{reduced Palm} distribution of $ T_{t,t+\epsilon}(\Phi)$ given a point at $x$. By Slivnyak's theorem, this reduced Palm distribution is the distribution of a PPP of intensity measure $\epsilon\Lambda$. Moreover, as $T_{t,t+\epsilon}(\Phi)$ is independent with $T_{t}\bigl(\mathcal{M}_{\infty}(\Phi)\bigr)$, its Palm distribution is also independent with the latter. Hence, by the refined Campbell formula,
\begin{align*}
& \mathbb{E}\left[\left|\Delta_{\infty,d,t,t+\epsilon}\cap A\right|\right]\\\noask
&=\,\epsilon\int_{A}\mathbb{E}_{x}^{!}\left[\prod_{y \in T_{t,t+\epsilon}(\Phi)}\textbf{1}_{C(x,y)=0}\right]\prod_{y \in T_{t}(\mathcal{M}_{\infty}(\Phi))}\textbf{1}_{C(x,y)=0}\,\Lambda(dx) \\\noask
&=\,\epsilon\int_{A}\mathbb{E}_{x}^{!}\left[\prod_{y \in  T_{t,t+\epsilon}(\Phi)}\bigl(1-h(x,y)\bigr)\right]\prod_{y \in T_{t}(\mathcal{M}_{\infty}(\Phi))}\bigl(1-h(x,y)\bigr)\,\Lambda(dx) \\\noask
&=\,\epsilon\int_{A}\exp\left \{-\epsilon\int_{\mathbb{R}^{2}}h(x,y)\,\Lambda(dy)\right \}\,\Lambda_{\infty,t}(dx). 
\end{align*}
 Thus, the intensity measure of $\Delta_{\infty,d,t,t+\epsilon}$ is 
\[\epsilon \exp\left \{-\epsilon\int_{\mathbb{R}^{2}}h(x,y)\,\Lambda(dy)\right \}\,\Lambda_{\infty,t}(dx).\]
We now apply Campbell's formula to $\Delta_{\infty,d,t,t+\epsilon}$,
\begin{align*}
&\mathbb{E}\hspace{-.1cm}\left[\sum_{x \in \Delta_{\infty,d,t,t+\epsilon}}\hspace{-.1cm}\bigl(1-v(x)\bigr)\right]\hspace{-.1cm}=\hspace{-.1cm} \epsilon\int_{\mathbb{R}^{2}}\bigl(1-v(x)\bigr)\exp\left \{-\epsilon\int_{\mathbb{R}^{2}}h(x,y)\,\Lambda(dy)\right \}
\Lambda_{\infty,t}(dx).
\end{align*}
As 
\begin{align*}
 \exp\left \{-\epsilon\int_{\mathbb{R}^{2}}h(x,y)\,\Lambda(dy)\right \}&\geq 1-\epsilon\int_{\mathbb{R}^{2}}h(x,y)\,\Lambda(dy) \geq 1-\epsilon\, \overline{H} ,
\end{align*}
we get
\begin{align}
 \nonumber &\epsilon\int_{\mathbb{R}^{2}}\bigl(1-v(x)\bigr)\,\Lambda_{\infty,t}(dx)\geq \mathbb{E}\left[\sum_{x \in T_{t,t+\epsilon}(\mathcal{M}_{\infty}(\Phi))}\bigl(1-v(x)\bigr)\right]\\\noask
\nonumber &\geq \epsilon\int_{\mathbb{R}^{2}}\bigl(1-v(x)\bigr)\,\Lambda_{\infty,t}(dx)-\epsilon^{2}\overline{H}\int_{\mathbb{R}^{2}}\bigl(1-v(x)\bigr)\,\Lambda_{\infty,t}(dx)\\\noask
 &\geq \epsilon\int_{\mathbb{R}^{2}}\bigl(1-v(x)\bigr)\,\Lambda_{\infty,t}(dx)-\epsilon^{2}\overline{H}\left(\int_{\mathbb{R}^{2}}\bigl(1-v(x)\bigr)\,\Lambda(dx)\right) \mbox{ a.s. } \label{E: thin lower bound}
 \end{align}
The conclusion then follows directly by putting together \eqref{E: thin upper bound} and \eqref{E: thin lower bound}.
\section{Proof of Theorem \ref{T:k Matern eq}}\label{A: Proof Theorem k Matern}
The proof that $g_{k,\Lambda}$ is continuous in $t$ is similar to the proof of Proposition \ref{P:f cont}.
 For the second part, following the same method as in the proof of Theorem \ref{T:CSMA eq}, we compute the conditional expectation
\begin{align*}
\mathbb{E}\left [\left .\left (\prod_{i=1}^{k}\left (\prod_{x \in T_{t,t+\epsilon}(\mathcal{Q}_{i}(\Phi))}v_{i}(x)\right )\right )\left (\prod_{x \in T_{t,t+\epsilon}(\mathcal{R}_{k}(\Phi))}v_{k+1}(x)\right )-1\,\right |\,  T_{t}(\Phi)\right ].
\end{align*}
To avoid cumbersome notation, we write $\mathbb{E}_{t}$ for the conditional expectation given $T_{t}(\Phi)$. 
By the same bounding technique as in Lemma \ref{L:limit 0}, we get 
\begin{align*}
&\left |1-\mathbb{E}_{t}\left [\left (\prod_{i=1}^{k}\left (\prod_{x \in T_{t,t+\epsilon}(\mathcal{Q}_{i}(\Phi))}v_{i}(x)\right )\right )\left (\prod_{x \in T_{t,t+\epsilon}(\mathcal{R}_{k}(\Phi))}v_{k}(x)\right )-1 \right ]+\right .\\\noask
&\left.\sum_{i=1}^{k}\mathbb{E}_{t}\left [\sum_{x \in T_{t,t+\epsilon}(\mathcal{Q}_{i}(\Phi))}\bigl(1-v_{i}(x)\bigr)\right ]+\mathbb{E}_{t}\left [\sum_{x \in T_{t,t+\epsilon}(\mathcal{R}_{k}(\Phi))}\bigl(1-v_{k+1}(x)\bigr)\right]\right |\\\noask
&\leq \sum_{i=2}^{\infty}\mathbb{E}_{t}\left [\sum_{x_{1},\ldots,x_{i} \textnormal{ mutually different in } T_{t,t+\epsilon}(\Phi)}\prod_{j=1}^{i}\bigl(1-v_{\min}(x_{j})\bigr)\right ]\\\noask
&=\sum_{i=2}^{\infty}\left (\int_{\mathbb{R}^{2}}\bigl(1-v_{\min}(x)\bigr)\,\epsilon\Lambda(dx)\right )^{i} \leq 2\epsilon^{2}\left (\int_{\mathbb{R}^{2}}\bigl(1-v_{\min}(x)\bigr)\,\Lambda(dx)\right ) \, \mbox{a.s.},
\end{align*}
where $v_{\min}(x)=\min_{i=1,\ldots,k+1}\bigl(v_{i}(x)\bigr)$ and $\epsilon$ small enough. In particular, one can take any $\epsilon$ smaller than $\left (2\int_{\mathbb{R}^{2}}\bigl(1-v_{\min}(x)\bigr)\,\Lambda(dx)\right )^{-1}$ to get
\begin{align}
\nonumber&\left |1-\mathbb{E}_{t}\left [\left (\prod_{i=1}^{k}\left (\prod_{x \in T_{t,t+\epsilon}(\mathcal{Q}_{i}(\Phi))}v_{i}(x)\right )\right )\left (\prod_{x \in T_{t,t+\epsilon}(\mathcal{R}_{k}(\Phi))}v_{k}(x)\right )-1 \right ]+\right .\\\noask
&\left.\sum_{i=1}^{k}\mathbb{E}_{t}\left [\sum_{x \in T_{t,t+\epsilon}(\mathcal{Q}_{i}(\Phi))}\bigl(1-v_{i}(x)\bigr)\right ]+\mathbb{E}_{t}\left [\sum_{x \in T_{t,t+\epsilon}(\mathcal{R}_{k}(\Phi))}\bigl(1-v_{k+1}(x)\bigr)\right]\right |\leq \epsilon
\end{align}

We now bound $\mathbb{E}_{t}\left [\displaystyle\sum_{x \in T_{t,t+\epsilon}(\mathcal{Q}_{i}(\Phi))}\bigl(1-v(x)\bigr)\right ]$ and $\mathbb{E}_{t}\left [\displaystyle\sum_{x \in T_{t,t+\epsilon}(\mathcal{R}_{i}(\Phi))}\bigl(1-v(x)\bigr)\right]$ with $v$ any function taking value in $[0,1]$ and satisfying \eqref{E: gen def cond}. We start with the case where $i$ is even, i.e. $i=2j$. First note that $\mathcal{Q}_{2j}(\Phi)=\mathcal{M}_{2j-2}(\Phi)\setminus \mathcal{M}_{2j}(\Phi)$ and $\mathcal{R}_{2j}(\Phi)=\mathcal{M}_{2j}(\Phi)\setminus \mathcal{M}_{2j-1}(\Phi)$ a.s.
Hence,
\begin{align*}
&\mathbb{E}_{t}\left [\sum_{x \in T_{t,t+\epsilon}(\mathcal{Q}_{2j}(\Phi))}\bigl(1-v(x)\bigr)\right ]\\\noask
&=\mathbb{E}_{t}\left [\sum_{x \in T_{t,t+\epsilon}(\mathcal{M}_{2j-2}(\Phi))}\bigl(1-v
(x)\bigr)\right ]-\mathbb{E}_{t}\left [\sum_{x \in T_{t,t+\epsilon}(\mathcal{M}_{2j}(\Phi))}\bigl(1-v(x)\bigr)\right ];\\\noask
&\mathbb{E}_{t}\left [\sum_{x \in T_{t,t+\epsilon}(\mathcal{R}_{2j}(\Phi))}\bigl(1-v(x)\bigr)\right ]\\\noask
&=\mathbb{E}_{t}\left [\sum_{x \in T_{t,t+\epsilon}(\mathcal{M}_{2j}(\Phi))}\bigl(1-v(x)\bigr)\right ]-\mathbb{E}_{t}\left [\sum_{x \in T_{t,t+\epsilon}(\mathcal{M}_{2j-1}(\Phi))}\bigl(1-v(x)\bigr)\right ] \mbox{a.s.}.
\end{align*}
By Proposition \ref{P:inf thin layer}, we have that $\Delta_{i,d,t,\epsilon}\subseteq T_{t,t+\epsilon}(\mathcal{M}_{i}(\Phi))\subseteq \Delta_{i,u,t,\epsilon}$ for every positive integer $i$ a.s. Hence,
\begin{align*}
&\mathbb{E}_{t}\left [\sum_{x \in \Delta_{2j-2,u,t,\epsilon}}\bigl(1-v(x)\bigr)\right ]-\mathbb{E}_{t}\left [\sum_{x \in \Delta_{2j,d,t,\epsilon}}\bigl(1-v(x)\bigr)\right ]\\\noask
&\geq \mathbb{E}_{t}\left [\sum_{x \in T_{t,t+\epsilon}(\mathcal{Q}_{2j}(\Phi))}\bigl(1-v(x)\bigr)\right ]\\\noask
&\geq \mathbb{E}_{t}\left [\sum_{x \in \Delta_{2j-2,d,t,\epsilon}}\bigl(1-v(x)\bigr)\right ]-\mathbb{E}_{t}\left [\sum_{x \in \Delta_{2j,u,t,\epsilon}}\bigl(1-v(x)\bigr)\right ];\\\noask
&\mathbb{E}_{t}\left [\sum_{x \in \Delta_{2j,u,t,\epsilon}}\bigl(1-v(x)\bigr)\right ]-\mathbb{E}_{t}\left [\sum_{x \in \Delta_{2j-1,d,t,\epsilon}}\bigl(1-v(x)\bigr)\right ]\\\noask
&\geq \mathbb{E}_{t}\left [\sum_{x \in T_{t,t+\epsilon}(\mathcal{R}_{2j}(\Phi))}\bigl(1-v(x)\bigr)\right ]\\\noask
&\geq\mathbb{E}_{t}\left [\sum_{x \in \Delta_{2j,d,t,\epsilon}}\bigl(1-v(x)\bigr)\right ]-\mathbb{E}_{t}\left [\sum_{x \in \Delta_{2j-1,u,t,\epsilon}}\bigl(1-v(x)\bigr)\right ] \mbox{ a.s. }
\end{align*}
We now compute $\mathbb{E}_{t}\left [\displaystyle\sum_{x \in \Delta_{i,u,t,\epsilon}}\bigl(1-v(x)\bigr)\right ]$ and bound $\mathbb{E}_{t}\left [\displaystyle\sum_{x \in \Delta_{i,d,t,\epsilon}}\bigl(1-v(x)\bigr)\right ]$. This is done by following the same bounding method in the proof of Lemma \ref{L:limit 0}, \begin{align*}
\mathbb{E}_{t}\left [\sum_{x \in \Delta_{i,u,t,\epsilon}}\bigl(1-v(x)\bigr)\right ]&=\epsilon\int_{\mathbb{R}^{2}}\bigl(1-v(x)\bigr)\,\Lambda_{i,t}(dx) ;\\\noask
\mathbb{E}_{t}\left [\sum_{x \in \Delta_{i,d,t,\epsilon}}\bigl(1-v(x)\bigr)\right ]&= \epsilon\int_{\mathbb{R}^{2}}\bigl(1-v(x)\bigr)e^{-\epsilon\overline{H}}\,\Lambda_{i,t}(dx)\\\noask
&\geq \epsilon\int_{\mathbb{R}^{2}}\bigl(1-v(x)\bigr)\left (1-\epsilon\overline{H}\right )\,\Lambda_{i,t}(dx)\\\noask
&\geq \epsilon\int_{\mathbb{R}^{2}}\bigl(1-v(x)\bigr)\,\Lambda_{i,t}(dx)-\epsilon^{2}\overline{H}\int_{\mathbb{R}^{2}}\bigl(1-v(x)\bigr)\,\Lambda(dx)\mbox{ a.s.},
\end{align*}
where $\Lambda_{i,t}(dx)=\prod_{y \in T_{t}(\mathcal{M}_{i-1}(\Phi))}\bigl(1-h(x,y)\bigr)\,\Lambda(dx)$. The last inequality comes from the fact that $\Lambda_{i,t}$ is smaller than $\Lambda$ almost everywhere. So,
\begin{align*}
&\epsilon\int_{\mathbb{R}^{2}}\bigl(1-v(x)\bigr)\,(\Lambda_{2j-2,t}-\Lambda_{2j,t})(dx)+\epsilon^{2}\overline{H}\int_{\mathbb{R}^{2}}\bigl(1-v(x)\bigr)\,\Lambda(dx)\\\noask
&\geq \mathbb{E}_{t}\left [\sum_{x \in T_{t,t+\epsilon}(\mathcal{Q}_{2j}(\Phi))}\bigl(1-v(x)\bigr)\right ]\\\noask
&\geq \epsilon\int_{\mathbb{R}^{2}}\bigl(1-v(x)\bigr)\,(\Lambda_{2j-2,t}-\Lambda_{2j,t})(dx)-\epsilon^{2}\overline{H}\int_{\mathbb{R}^{2}}\bigl(1-v(x)\bigr)\,\Lambda(dx);\\\noask
&\epsilon\int_{\mathbb{R}^{2}}\bigl(1-v(x)\bigr)\,(\Lambda_{2j,t}-\Lambda_{2j-1,t})(dx)+\epsilon^{2}\overline{H}\int_{\mathbb{R}^{2}}\bigl(1-v(x)\bigr)\,\Lambda(dx)\\\noask
&\geq\mathbb{E}\left [\sum_{x \in T_{t,t+\epsilon}(\mathcal{R}_{2j}(\Phi))}\bigl(1-v(x)\bigr)\right ]\\\noask
&\geq \epsilon\int_{\mathbb{R}^{2}}\bigl(1-v(x)\bigr)\,(\Lambda_{2j,t}-\Lambda_{2j-1,t})(dx)-\epsilon^{2}\overline{H}\int_{\mathbb{R}^{2}}\bigl(1-v(x)\bigr)\,\Lambda(dx).
\end{align*} 
Doing similarly for the second case where $i=2j-1$, we get
\begin{align*}
&\epsilon\int_{\mathbb{R}^{2}}\bigl(1-v(x)\bigr)\,(\Lambda_{2j-1,t}-\Lambda_{2j-3,t})(x)+\epsilon^{2}\overline{H}\int_{\mathbb{R}^{2}}\bigl(1-v(x)\bigr)\,\Lambda(dx)\\\noask
&\geq\mathbb{E}_{t}\left [\sum_(x \in T_{t,t+\epsilon}{\mathcal{Q}_{2j-1}(\Phi))}\bigl(1-v(x)\bigr)\right ]\\\noask
&\epsilon\int_{\mathbb{R}^{2}}\bigl(1-v(x)\bigr)\,(\Lambda_{2j
-1,t}-\Lambda_{2j-3,t})(x)-\epsilon^{2}\overline{H}\int_{\mathbb{R}^{2}}\bigl(1-v(x)\bigr)\,\Lambda(dx)\\\noask
&\epsilon\int_{\mathbb{R}^{2}}\bigl(1-v(x)\bigr)\,(\Lambda_{2j-2,t}-\Lambda_{2j-1,t})(x)+\epsilon^{2}\overline{H}\int_{\mathbb{R}^{2}}\bigl(1-v(x)\bigr)\,\Lambda(dx)\\\noask
&\geq\mathbb{E}_{t}\left [\sum_{x \in T_{t,t+\epsilon}(\mathcal{R}_{2j-1}(\Phi))}\bigl(1-v(x)\bigr)\right ]\\\noask
&\epsilon\int_{\mathbb{R}^{2}}\bigl(1-v(x)\bigr)\,(\Lambda_{2j-2,t}-\Lambda_{2j-1,t})(x)-\epsilon^{2}\overline{H}\int_{\mathbb{R}^{2}}\bigl(1-v(x)\bigr)\,\Lambda(dx).
\end{align*}
Note that in the above formula, $\Lambda_{-1}(dx)=0$ by convention.
So,
\begin{align}
\nonumber&-\lim_{\epsilon\rightarrow 0}\epsilon^{-1}\mathbb{E}_{t}\left [\left (\prod_{i=1}^{2k}\left (\prod_{x \in T_{t,t+\epsilon}(\mathcal{Q}_{i}(\Phi))}v_{i}(x)\right )\right )\left (\prod_{x \in T_{t,t+\epsilon}(\mathcal{R}_{2k}(\Phi))}v_{2k+1}(x)\right )-1\,\right ]\\
\nonumber&=\int_{\mathbb{R}^{2}}\bigl(1-v_{2k+
1}(x)\bigr)\,(\Lambda_{2k,t}-\Lambda_{2k-1,t})(dx)+ \sum_{j=1}^{k}\bigl(1-v_{2j-1}(x)\bigr)\,(\Lambda_{2j-1,t}-\Lambda_{2j-3,t})(dx)\\
\nonumber&\hspace{.5 cm}+\sum_{j=1}^{k}\int_{\mathbb{R}^{2}}\bigl(1-v_{2j}(x)\bigr)\,(\Lambda_{2j-2,t}-\Lambda_{2j,t})(dx)\\\noask
\nonumber&=\sum_{j=1}^{k}\int_{\mathbb{R}^{2}}\bigl(v_{2j+1}(x)-v_{2j-1}(x)\bigr)\,\Lambda_{2j-1,t}(dx)+\sum_{j=1}^{k}\int_{\mathbb{R}^{2}}\bigl(v_{2j}(x)-v_{2j+2}(x)\bigr)\,\Lambda_{2j,t}(dx)\\
&=\sum_{i=1}^{2k}\bigl(w_{i}(x)-v_{i}(x)\bigr)\,\Lambda_{i,t}(dx);\label{E: lim 2k}
\end{align}
and
\begin{align}
\nonumber&-\lim_{\epsilon\rightarrow 0}\epsilon^{-1}\mathbb{E}_{t}\left [\left (\prod_{i=1}^{2k-1}\left (\prod_{x \in T_{t,t+\epsilon}(\mathcal{Q}_{i}(\Phi))}v_{i}(x)\right )\right )\left (\prod_{x \in T_{t,t+\epsilon}(\mathcal{R}_{2k-1}(\Phi))}v_{2k}(x)\right )-1\,\right ]\\
\nonumber&=\int_{\mathbb{R}^{2}}\bigl(1-v_{2k}(x)\bigr)\,(\Lambda_{2k-2,t}-\Lambda_{2k-1,t})(dx)+ \sum_{j=1}^{k}\int_{\mathbb{R}^{2}}\bigl(1-v_{2j-1}(x)\bigr)\,(\Lambda_{2j-1,t}-\Lambda_{2j-3,t})(dx)\\
\nonumber&\hspace{.5 cm}+\sum_{j=1}^{k-1}\int_{\mathbb{R}^{2}}\bigl(1-v_{2j}(x)\bigr)\,(\Lambda_{2j-2,t}-\Lambda_{2j,t})(dx)\\\noask
\nonumber&=\sum_{j=1}^{k-1}\int_{\mathbb{R}^{2}}\bigl(v_{2j+1}(x)-v_{2j-1}(x)\bigr)\,\Lambda_{2j-1,t}(dx)+\sum_{j=1}^{k-1}\int_{\mathbb{R}^{2}}\bigl(v_{2j}(x)-v_{2j+2}(x)\bigr)\,\Lambda_{2j,t}(dx)\\\noask
&=\sum_{i=1}^{2k+1}\bigl(w_{i}(x)-v_{i}(x)\bigr)\,\Lambda_{i,t}(dx) \mbox{ a.s.}\label{E: lim 2k+1}
\end{align}
By noting that 
\begin{align*}
\mathcal{M}_{i}(\Phi)&=\left (\bigcup_{k\geq j\geq 1, i \in I_{j} }\mathcal{Q}_{j}(\Phi)\right )\cup \mathcal{R}_{k}(\Phi)
\end{align*}
if $i \in J_{k}$ and
\begin{align*}
\mathcal{M}_{i}(\Phi)&=\left (\bigcup_{k\geq j\geq 1, i \in I_{j} }\mathcal{Q}_{j}(\Phi)\right )
\end{align*}
otherwise, we have
\begin{align}
\nonumber&\mathbb{E}\left [\left (\prod_{j=1}^{k}\left (\prod_{y \in T_{t}(\mathcal{Q}_{j}(\Phi))}v_{j}(y)\right )\right )\left (\prod_{y \in T_{t}(\mathcal{R}_{k}(\Phi))}v_{k+1}(y)\right )\Lambda_{i,t}(dx)\right ]\\\noask
\nonumber&=\mathbb{E}\left  [\left (\prod_{j=1}^{k}\left (\prod_{y \in T_{t}(\mathcal{Q}_{j}(\Phi))}v_{j}(y)\right )\right )\left (\prod_{y \in T_{t}(\mathcal{R}_{k}(\Phi))}v_{k+1}(y)\right )\right .\\
\nonumber&\hspace{.5cm}\left .\left (\prod_{y \in T_{t}(\mathcal{M}_{i-1}(\Phi))}\bigl(1-h(x,y)\bigr)\right )\Lambda(dx)\right  ]\\\noask
\nonumber&= \mathbb{E}\left  [\left (\prod_{j=1}^{k}\left (\prod_{y \in T_{t}(\mathcal{Q}_{j}(\Phi))}v_{j}(y)\right )\right )\left (\prod_{y \in T_{t}(\mathcal{R}_{k}(\Phi))}v_{k+1}(y)\right )\right .\\
\nonumber&\hspace{.5cm}\left .\left (\prod_{k \geq j \geq 1, i-1 \in I_{j}}\left (\prod_{y \in T_{t}(\mathcal{Q}_{j}(\Phi))}\bigl(1-h(x,y)\bigr)\right )\prod_{y \in T_{t}(\mathcal{R}_{k}(\Phi))}\bigl(1-h(x,y)\bigr)\right )\right  ]\Lambda(dx)\\\noask
\nonumber&=\mathbb{E}\left  [\left (\prod_{j=1}^{k}\left (\prod_{y \in T_{t}(\mathcal{Q}_{i}(\Phi))}\mathfrak{u}_{j,i,k}(y,x)\right )\right )\left (\prod_{y \in T_{t}(\mathcal{R}_{k}(\Phi))}\mathfrak{u}_{k+1,i,k}(y,x)\right )\right  ]\Lambda(dx)\\\noask
&=g_{k,\Lambda}(t,\textbf{H}_{i,k}(\textbf{v},x))
\end{align}
if $i-1 \in J_{k}$ and 
\begin{align}
\nonumber&\mathbb{E}\left [\left (\prod_{j=1}^{k}\left (\prod_{y \in T_{t}(\mathcal{Q}_{j}(\Phi))}v_{j}(y)\right )\right )\left (\prod_{y \in T_{t}(\mathcal{R}_{k}(\Phi))}v_{k+1}(y)\right )\Lambda_{i,t}(dx)\right ]\\\noask
\nonumber&=\mathbb{E}\left  [\left (\prod_{j=1}^{k}\left (\prod_{y \in T_{t}(\mathcal{Q}_{j}(\Phi))}v_{j}(y)\right )\right )\left (\prod_{y \in T_{t}(\mathcal{R}_{k}(\Phi))}v_{k+1}(y)\right )\right .\\
\nonumber&\hspace{.5cm}\left .\left (\prod_{y \in T_{t}(\mathcal{M}_{i-1}(\Phi))}\bigl(1-h(x,y)\bigr)\right )\Lambda(dx)\right  ]\\\noask
\nonumber&= \mathbb{E}\left  [\left (\prod_{j=1}^{k}\left (\prod_{y \in T_{t}(\mathcal{Q}_{j}(\Phi))}v_{j}(y)\right )\right )\left (\prod_{y \in T_{t}(\mathcal{R}_{k}(\Phi))}v_{k+1}(y)\right )\right .\\
\nonumber&\hspace{.5cm}\left .\left (\prod_{k \geq j \geq 1, i-1 \in I_{j}}\left (\prod_{y \in T_{t}(\mathcal{Q}_{j}(\Phi))}\bigl(1-h(x,y)\bigr)\right )\right )\right  ]\Lambda(dx)\\\noask
\nonumber&=\mathbb{E}\left  [\left (\prod_{j=1}^{k}\left (\prod_{y \in T_{t}(\mathcal{Q}_{i}(\Phi))}\mathfrak{u}_{j,i,k}(y,x)\right )\right )\left (\prod_{y \in T_{t}(\mathcal{R}_{k}(\Phi))}\mathfrak{u}_{k+1,i,k}(y,x)\right )\right  ]\Lambda(dx)\\\noask
&=g_{k,\Lambda}(t,\textbf{H}_{i,k}(\textbf{v},x))
\end{align}
otherwise.

 We conclude by substituting the above equalities to \eqref{E: lim 2k} and \eqref{E: lim 2k+1}. 
 
\section{Proof of Proposition \ref{P:Palm gen infty Matern}}\label{A: Proof Proposition Palm}
 As there is no ambiguity, we drop the $\infty$ and $\Lambda$ subscripts in this proof.   By Proposition \ref{P:Palm non Palm}, 
\begin{align*}
 \frac{d}{ds}f(t,ve^{-s\textbf{1}_{B}})\Bigr|_{s=0}=-\int_{B}v(y)f_{y}(t,v)m(t,y)\Lambda(dy).
\end{align*}
By Theorem \ref{T:CSMA eq}, 
\begin{align*}
& f(t,ve^{-s\textbf{1}_{B}})=1-\int_{0}^{t}\int_{\mathbb{R}^{2}}f\bigl(\tau,H(ve^{-s\textbf{1}_{B}},x)\bigr)\left(1-v(x)e^{-s\textbf{1}_{B}(x)}\right)\,\Lambda(dx)\,d\tau.
\end{align*}
By the same argument as in the proof of Proposition \ref{P: First moment measure}, we can change the order of integration and derivation in the above formula,
\begin{align*}
 &\int_{B}v(y)f_{y}(t,v)m(t,y)\,\Lambda(dx)\\\noask
&=\frac{d}{ds}\left(\int_{0}^{t}\int_{\mathbb{R}^{2}}f\bigl(\tau,H(ve^{-s\textbf{1}_{B}},x)\bigr)\left(1-v(x)e^{-s\textbf{1}_{B}(x)}\right)\,\Lambda(dx)\,d\tau\right)\Biggr|_{s=0}\\\noask
&=\int_{0}^{t}\int_{\mathbb{R}^{2}}\frac{d}{ds}\Bigl(f\bigl(\tau,H(ve^{-s\textbf{1}_{B}},x)\bigr)\left(1-v(x)e^{-s\textbf{1}_{B}(x)}\right)\Bigr)\Biggr|_{s=0}\,\Lambda(dx)\,d\tau.
\end{align*}
Moreover,
\begin{align*}
 &\frac{d}{ds}\Bigl(f\bigl(\tau,H(ve^{-s\textbf{1}_{B}},x)\bigr)\left(1-v(x)e^{-s\textbf{1}_{B}(x)}\right)\Bigr)\Bigr|_{s=0}\\\noask
 &=\frac{d}{ds}f\bigl(\tau,H(ve^{-s\textbf{1}_{B}},x)\bigr)\Bigr|_{s=0}\bigl(1-v(x)\bigr)+\textbf{1}_{B}(x)v(x)f(\tau,H(v,x)).
\end{align*}
By Proposition \ref{P:Palm non Palm}, 
\begin{align*}
  \frac{d}{ds}f\bigl(\tau,H(ve^{-s\textbf{1}_{B}},x)\bigr)\Bigr|_{s=0}&=\frac{d}{ds}f\bigl(\tau,H(v,x)e^{-s\textbf{1}_{B}}\bigr)\Bigr|_{s=0}\\\noask
&=-\int_{B}f_{y}\bigl(\tau,H(v,x)\bigr)v(y)\bigl(1-h(x,y)\bigr)m(\tau,y)\,\Lambda(dy)
\end{align*}
for every $x$. Hence,
\begin{align*}
 &\int_{B}v(y)f_{y}(t,v)m(t,y)\Lambda(dx)\\\noask
&=-\int_{0}^{t}\int_{\mathbb{R}^{2}}\int_{B}f_{y}\bigl(\tau,H(v,x)\bigr)v(y)\bigl(1-h(x,y)\bigr)m(\tau,y)\,\Lambda(dy)\,\Lambda(dx)\,d\tau\\
&\hspace{.5cm}+\int_{0}^{t}\int_{B}f\bigl(\tau,H(v,x)\bigr)v(x)\,\Lambda(dx)\,d\tau\\\noask
&=-\int_{0}^{t}\int_{\mathbb{R}^{2}}\int_{B}f_{y}\bigl(\tau,H(v,x)\bigr)v(y)\bigl(1-h(x,y)\bigr)m(\tau,y)\,\Lambda(dy)\,\Lambda(dx)\,d\tau\\
&\hspace{.5cm}+\int_{0}^{t}\int_{B}f\bigl(\tau,H(v,y)\bigr)v(y)\,\Lambda(dy)\,d\tau\\
&=-\int_{B}\int_{0}^{t}\int_{\mathbb{R}^{2}}f_{y}\bigl(\tau,H(v,x)\bigr)v(y)\bigl(1-h(x,y)\bigr)m(\tau,y)\,\Lambda(dx)\,d\tau\,\Lambda(dy)\\
&\hspace{.5cm}+\int_{B}\int_{0}^{t}f\bigl(\tau,H(v,y)\bigr)v(y)\,d\tau
\,\Lambda(dy).
\end{align*}
As this is true for any bounded measurable set $B$ such that $v(x)>0$ for all $x$ in $B$, we must have
\begin{align*}
 f_{y}(t,v)&=\frac{\int_{0}^{t}f(\tau,H(v,y))d\tau}{m(t,y)}\\
&\hspace{.5cm}-\int_{0}^{t}\int_{\mathbb{R}^{2}}f_{y}\bigl(\tau,H(v,x)\bigr)\bigl(1-v(x)\bigr)\bigl(1-h(x,y)\bigr)\frac{m(\tau,y)}{m(t,y)}\Lambda(dx)d\tau
\end{align*}
for $\Lambda$- almost every $y$ such that $v(y)>0$.
For $y$ such that $v(y)=0$, by Proposition \ref{P:reduced Palm non Palm}, for any bounded measurable set $B$ not in the support of $v$,
\begin{align*}
\frac{d}{ds}f(t,v+s\textbf{1}_{B})\Bigr|_{s=0}=\int_{B}f_{y}(t,v)m(t,y)\Lambda(dy). 
\end{align*}
Again, by Theorem \ref{T:CSMA eq},
\begin{align*}
& f(t,v+s\textbf{1}_{B})\\\noask
&=1-\int_{0}^{t}\int_{\mathbb{R}^{2}}f\bigl(\tau,H(v+s\textbf{1}_{B},x)\bigr)\bigl(1-v(x)-s\textbf{1}_{B}(x)\bigr)\,\Lambda(dx)\,d\tau.
\end{align*}
By using the same arguments as above and by noting that
\begin{align*}
  \frac{d}{ds}f\bigl(\tau,H(v+s\textbf{1}_{B},x)\bigr)\Bigr|_{s=0}&=\frac{d}{ds}f\bigl(\tau,H(v,x)+s\textbf{1}_{B}(1-h(.,x))\bigr)\Bigr|_{s=0}\\\noask
&=\int_{B}f_{y}\bigl(\tau,H(v,x)\bigr)\bigl(1-h(x,y)\bigr)m(\tau,y)\,\Lambda(dy),
\end{align*}
we have
\begin{align*}
 f_{y}(t,v)&=\frac{\int_{0}^{t}f\bigl(\tau,H(v,y)\bigr)d\tau}{m(t,y)}\\
&\hspace{.5cm}-\int_{0}^{t}\int_{\mathbb{R}^{2}}f_{y}\bigl(\tau,H(v,x)\bigr)\bigl(1-v(x)\bigr)\bigl(1-h(x,y)\bigr)\frac{m(\tau,y)}{m(t,y)}\Lambda(dx)d\tau
\end{align*}
for $\Lambda$- almost every $x$ such that $v(x)=0$.

\end{document}